\documentclass[12pt]{amsart}

\usepackage{DKstyle}

\newcommand{\Mod}[1]{\ (\mathrm{mod}\ #1)}
\newcommand{\vertiii}[1]{{\left\vert\kern-0.25ex\left\vert\kern-0.25ex\left\vert #1
    \right\vert\kern-0.25ex\right\vert\kern-0.25ex\right\vert}}
\theoremstyle{plain}

\newtheorem*{theorem*}{Theorem}
\newcommand*{\medcap}{\mathbin{\scalebox{1.5}{\ensuremath{\cap}}}}
\makeatletter
\newcommand*{\rom}[1]{\expandafter\@slowromancap\romannumeral #1@}
\makeatother

\usepackage{bm}
\usepackage{caption}
\usepackage[labelfont=rm]{subcaption}

\usepackage[pagebackref=true, colorlinks]{hyperref}

\hypersetup{pdffitwindow=true,linkcolor=blue,citecolor=blue,urlcolor=blue,menucolor=blue}

\usepackage{comment}


\date{\today}%
\dedicatory{}%
\commby{}%

\title{Effective density for inhomogeneous quadratic forms I: generic forms and fixed shifts}
\author{Anish Ghosh}
\address{School of Mathematics, Tata Institute of Fundamental Research, Homi Bhabha Road, Colaba, Mumbai 400005, India}
\email{ghosh@math.tifr.res.in}
\author{Dubi Kelmer}
\address{Department of Mathematics, Boston College, Chestnut Hill MA 02467-3806, USA}
\email{kelmer@bc.edu}
\author{Shucheng Yu}
\address{Department of Mathematics, Technion, Haifa, Israel}
\email{yushucheng@campus.technion.ac.il}
\thanks{AG was supported by a Government of India, Department of Science and Technology, Swarnajayanti fellowship, a CEFIPRA grant, a MATRICS grant and a grant from the Infosys Foundation. AG acknowledges support from the Department of Atomic Energy, Government of India, under project $12-R\&D-TFR-5.01-0500$. DK and SY were partially supported by NSF CAREER grant DMS-1651563. SY acknowledge the support of ISF grant number  871/17. SY acknowledges that this project has received funding from the European Research Council (ERC) under the European Union's Horizon 2020 research and innovation program (grant agreement No.\ 754475).}
\begin{document}
\begin{abstract}
We establish effective versions of Oppenheim's conjecture for generic inhomogeneous quadratic forms. We prove such results for fixed shift vectors and generic quadratic forms. When the shift is rational we prove a counting result which implies the optimal density for values of generic inhomogeneous forms. We also obtain a similar density result for fixed irrational shifts satisfying an explicit Diophantine condition. The main technical tool is  a formula for the second moment of Siegel transforms on certain congruence quotients of $\SL_n(\R)$ which we believe to be of independent interest. In a sequel, we use different techniques to treat the companion problem concerning generic shifts and fixed quadratic forms.
\end{abstract}
\maketitle

\tableofcontents
\section{Introduction}
Let $Q$ be a quadratic form on $\R^{n}$ and let $\bm{\alpha}$ be a vector in $\R^{n}$. Define the inhomogeneous quadratic form $Q_{\bm{\alpha}}$ by
$$Q_{\bm{\alpha}}(\bm{v})=Q(\bm{v}+\bm{\alpha})\ \textrm{for any}\ \bm{v}\in \R^{n},$$
where we think of $Q_{\bm\alpha}$ as a \textit{shift} by ${\bm \alpha}$ of the homogenous form $Q$.
The inhomogeneous form $Q_{\bm{\alpha}}$ is said to be \textit{indefinite} if $Q$ is indefinite and \textit{non-degenerate} if $Q$ is non-degenerate. Finally, $Q_{\bm{\alpha}}$ is said to be \textit{irrational} if either $Q$ is an irrational quadratic form, i.e. not proportional to a quadratic form with integer coefficients, or $\bm{\alpha}$ is an irrational vector. 

In \cite{MargulisMohammadi11},  Margulis and Mohammadi studied the counting function 
\begin{equation}\label{equ:countingfunction}
\cN_{Q_{\bm\alpha},I}(t)=\#\{\bm{v}\in\Z^n\ |\  Q_{\bm \alpha}(\bm{v})\in I,\ \|\bm{v}\|\leq t\},
\end{equation}
where $I\subseteq \R$ is an interval and $\|\cdot\|$ is the Euclidean norm on $\R^n$ \footnote{Their result is more general and deals with dilations of arbitrary convex sets.}.
They showed that for any indefinite, irrational and non-degenerate inhomogeneous form $Q_{\bm \alpha}$ in $n \geq 3$ variables there is $c_Q>0$ such that 
$$\liminf_{t\to \infty}\frac{\cN_{Q_{\bm \alpha},I}(t)}{t^{n-2}}\geq c_Q |I|,$$
while for $n\geq 5$ the limit exists and equals $c_Q|I|$. This generalizes the results of Dani--Margulis \cite{DaniMargulis1993} (a lower bound for the liminf) and Eskin--Margulis--Mozes \cite{EskinMargulisMozes1998
}  (an upper bound for the limsup) in the homogeneous case $\bm{\alpha}=\bm{0}$. It also generalizes results of Sarnak \cite{Sarnak1997} (for the homogeneous case), and Marklof \cite{Marklof2002,Marklof2003} (for the inhomogeneous case) who considered an important special case related to the pair correlation of values of a positive definite form. In particular, one can deduce from this an 
analogue of the famous Oppenheim conjecture 
for inhomogeneous forms, stating that for any indefinite, irrational,  non-degenerate inhomogeneous form $Q_{\bm \alpha}$ in $n \geq 3$ variables $Q_{\bm{\alpha}}(\mathbb{Z}^{n})$ is dense in $\mathbb{R}$ (see also \cite{BandiGhosh} for a self-contained proof).


In this paper, we are concerned with questions of \emph{effectivity}, namely, for a given $Q_{\bm{\alpha}}, \xi \in \R$ and $t \geq 1$ large:
\begin{enumerate}
\item How small can $|Q_{\bm{\alpha}}(\bm{v}) - \xi|$ get for  $\bm{v} \in \Z^{n}$ with $\|\bm{v}\|\leq t$ bounded? 
\item Is it possible to obtain an effective estimate  (i.e., with power saving for the remainder)  for the counting function $\cN_{Q_{\bm\alpha},I}(t)$?
\end{enumerate}
We note that the two problems are closely related and that a good enough answer for the second question will also give an answer to the first one.
\subsection{Homogeneous forms}
Before turning to our results we discuss what is known for homogeneous forms, a problem that received a lot of attention in recent years. When the number of variables $n\geq 5$ is large, it is possible to use analytic methods to get effective results depending on Diophantine properties of the form (see  \cite{GotzeMargulis2010} for the best results in this setting). For $n\leq 4$ the analytic methods break down and one has to rely on a dynamical approach. In this case, using an effective version of the results of Dani and Margulis \cite{DaniMargulis1993} on unipotent flows,  Lindenstrauss and Margulis  \cite{LindenstraussMargulis2014} answered the first question with a bound that is logarithmic in $t$ (under some mild explicit assumptions on the coefficients of $Q$). Their result is remarkable because it applies to an explicit, large class of quadratic forms. However, as there can be considerable variation in the Diophantine properties of quadratic forms, it is plausible to expect much better (polynomial) bounds for generic forms.
Indeed, in \cite{GhoshGorodnikNevo2018}, the problem was considered for generic forms and the authors gave a heuristic argument predicting that for generic forms one should expect that for all $\kappa<n-2$ the system of inequalities
\begin{equation}\label{equ:homoge}
|Q(\bm{v})-\xi|<t^{-\kappa},\; \|\bm{v}\|\leq t
\end{equation}
has integer solutions $\bm{v}\in \Z^n$ for all sufficiently large $t$. 

To make the notion of generic forms more precise, we note that any quadratic form of signature $(p_1,p_2)$ with $p_1+p_2=n$ is of the form $\lambda Q(\bm{v}g)$ with $Q$ some fixed unit determinant form of signature $(p_1,p_2)$, $g\in \SL_n(\R)$ and $\lambda>0$. We say that a property holds for almost all forms of signature $(p_1,p_2)$ if for any $\lambda>0$ it holds for $\lambda Q(\bm{v} g)$ for almost all $g\in \SL_n(\R)$, and that it holds for almost all indefinite forms in $n$ variables if it holds for almost all forms of signature $(p_1,p_2)$ for any $p_1, p_2\geq 1$ with $p_1+p_2=n$. With this notion, using an effective mean ergodic theorem for semisimple groups, Ghosh, Gorodnik and Nevo \cite{GhoshGorodnikNevo2018} showed that there is some positive $\kappa_0$ (depending on $n$) such that for any $\kappa<\kappa_0$ and any $\xi\in \R$ for almost all indefinite forms $Q$ in $n\geq 3$ variables \eqref{equ:homoge} has integer solutions for all sufficiently large $t$. For the special case of $n=3$ their method gives $\kappa_0=1$ as predicted (see also \cite{GhoshKelmer17}), but for $n>3$ the value of $\kappa_0$ is strictly smaller. The method also gives the predicted heuristic value for several other examples of polynomial maps on homogeneous varieties of semisimple groups.
By using a completely different method, relying on Rogers' formula \cite{Rogers1955} for the second moment of Siegel transforms, in \cite{AthreyaMargulis2018} Athreya and Margulis obtained the optimal rate $\kappa_0=n-2$ for any $n\geq 3$ in the special case of $\xi=0$. Moreover, they showed that for almost all indefinite forms in $n\geq 3$ variables, and for any fixed interval $I$, the counting function $\cN_{Q,I}(T)$ satisfies an asymptotic formula with a power saving estimate for the remainder.  Subsequently, Kelmer and Yu further developed this technique in \cite{KelmerYu2018b} allowing the interval $I=I_t$ to be shrinking intervals of of size $t^{-\kappa}$ for any $\kappa<n-2$, thus confirming the prediction made in \cite{GhoshGorodnikNevo2018}. The method used in \cite{KelmerYu2018b} is rather soft and only uses the second moment formula together with some explicit volume estimates, in particular it could be applied in any setting where these two ingredients are available.
\begin{rem}
We note that for these results, one first fixes the target point $\xi$ (or the shrinking interval $I_t$) and then picks a full measure set of forms with values approaching $\xi$. One can also ask what is the best rate at which a generic form can approximate all targets simultaneously. This type of problem was studied by  Bourgain \cite{Bourgain2016}, who proved a certain uniform effective result for generic diagonal forms in three variables using analytic methods. Subsequently, Ghosh and Kelmer \cite{GhoshKelmer2018} used ergodic methods to get similar results for generic ternary forms, and finally in \cite{KelmerYu2018b} Kelmer and Yu proved an analogous counting result in the uniform setting for generic forms in any number of variables.
\end{rem}


\subsection{Inhomogeneous forms}
Turning to inhomogeneous forms, as pointed out in \cite{AthreyaMargulis2018}, by using an affine analogue of Rogers' formula obtained  in \cite[Appendix B]{ElBazMarklofVinogradov2015} (see also \cite[Lemma 4]{Athreya2015}), one can recover the main results in \cite{AthreyaMargulis2018,KelmerYu2018b} with $Q$ replaced by $Q_{\bm{\alpha}}$, for almost all  indefinite forms $Q$ in $n$ variables and almost all $\bm{\alpha}\in \R^n$.
Since all the arguments are identical to the ones in \cite{KelmerYu2018b} (with the obvious modifications, replacing the space of lattices $\SL_n(\Z)\bk \SL_n(\R)$ by the space of affine lattices $\rm{ASL}_n(\Z)\bk \rm{ASL}_n(\R)$) we will not give the details here, but suffice to point out that it implies the following result answering both questions in this setting:
\begin{Thm}
Let $n\geq 3$ and let $0\leq \kappa<n-2$. Let $\{I_t\}_{t> 0}$ be a decreasing family of bounded measurable subsets of $\R$ with measures $|I_t|=ct^{-\kappa}$ for some $c>0$. Then there is $\nu>0$ such that for almost every non-degenerate indefinite quadratic form $Q$ in $n$ variables there is a constant $c_Q$ such that for almost every $\bm\alpha\in \R^n$
$$\cN_{Q_{\bm{\alpha},I_t}}(t)=c_Q |I_t|t^{n-2}+O_{Q,\bm{\alpha}}(t^{n-2-\kappa-\nu}).$$
In particular, for any $\kappa<n-2,$ for almost every $\bm\alpha\in\R^n,$ and for almost every quadratic form $Q$ as above, the system of inequalities 
\begin{equation} \label{equ:maininequ} 
|Q_{\bm{\alpha}}(\bm{v})-\xi|<t^{-\kappa},\; \|\bm{v}\|\leq t
\end{equation}
has integer solutions for all sufficiently large $t$.
\end{Thm}
\begin{rmk}\label{rmk:scale}
We note that, more precisely, following the arguments in \cite{KelmerYu2018b} we only get the above counting result for generic non-degenerate indefinite quadratic forms with unit determinant. However, it can be easily generalized to quadratic forms of all determinants by noting that $(\lambda Q_{\bm{\alpha}})^{-1}(I)=Q_{\bm{\alpha}}^{-1}(\lambda^{-1} I)$ and setting $c_{\lambda Q}=\lambda^{-1}c_Q$.
\end{rmk}
The main result of this paper is to address the significantly more subtle problem of obtaining similar estimates for a fixed shift $\bm\alpha$ (perhaps satisfying some Diophantine assumptions) that will hold for $Q_{\bm\alpha},$ for almost all forms $Q$. In a companion paper  \cite{GKY2020}, we address the companion problem of obtaining estimates that hold for a fixed form and for almost all shifts.  We should also mention the recent result of Str\"{o}mbergsson and Vishe \cite{StrombergssonVishe2018} who gave an effective estimate for $\cN_{Q_{\bm\alpha},I}(t)$ for the special case when $Q$ is the fixed standard form of signature $(2,2)$ and $\bm\alpha$ is a fixed shift satisfying some explicit Diophantine conditions.

\subsection{Rational shifts}
Our strongest result in this setting is when the shift $\bm{\alpha}$ is rational.
Here we can use a similar approach to \cite{KelmerYu2018b} to prove the following counting result.
\begin{Thm}\label{thm:rationalshift}
Let $n\geq 3$ and let $0\leq\kappa<n-2$. Let $\bm{\alpha}\in \Q^n$ be a fixed rational vector. Let $\{I_t\}_{t> 0}$ be a decreasing family of bounded measurable subsets of $\R$ with measures $|I_t|=ct^{-\kappa}$ for some $c>0$. Then there is $\nu>0$ such that for almost every non-degenerate indefinite quadratic form $Q$ in $n$ variables, there is a constant $c_Q>0$ such that 
\begin{equation}\label{equ:countingresult}
\cN_{Q_{\bm{\alpha},I_t}}(t)=c_Q |I_t|t^{n-2}+O_{Q,\bm{\alpha}}(t^{n-2-\kappa-\nu}).
\end{equation}
In particular, for any $\kappa<n-2,$ for every $\bm\alpha\in\Q^n,$ and for almost every quadratic form $Q$ as above, the system of inequalities \eqref{equ:maininequ} has integer solutions for all sufficiently large $t$.
\end{Thm}

Note that the set of rational vectors is a measure zero subset in $\R^n$ so the result on generic shifts does not say anything about rational shifts. Nevertheless, our result for a fixed rational shift gives the same (optimal) rate which holds generically. The reason we are able to obtain such strong results for a rational shift $\bm\alpha=\frac{\bm p}{q}$ is the following observation:
The values of the shifted form at integer points 
\begin{equation}\label{equ:obs}
Q_{\bm\alpha}(\bm v)=Q(\bm v+\frac{\bm p}{q})=q^{-2}Q(q\bm v+\bm{p}),
\end{equation}
are just a scaling of the values of the homogeneous form $Q$ evaluated on integer points satisfying a congruence condition modulo $q$. This observation also gives another interpretation of our counting result in terms of counting integer solutions to \eqref{equ:homoge} with an extra congruence condition. Explicitly, we show that, for a generic form $Q$, the integer solutions to \eqref{equ:homoge} are evenly distributed in each congruence class in $(\Z/q\Z)^n$.
\begin{Cor}\label{c:cong}
Let $n\geq 3$ and let $0\leq\kappa<n-2$ and let $q\in \N$. Let $\{I_t\}_{t> 0}$ be a decreasing family of bounded measurable subsets of $\R$ with measures $|I_t|=ct^{-\kappa}$ for some $c>0$. Then there is $\nu>0$ such that for almost every non-degenerate indefinite quadratic form $Q$ in $n$ variables, there is a constant $c_Q>0$ such that for all $\bm p\in (\Z/q\Z)^n$
$$\#\{\bm{w}\in \Z^n\ |\ Q(\bm{w})\in I_t, \bm{w}\equiv \bm{p} \Mod{q}, \|\bm{w}\|<t\}=\frac{c_Q|I_t|t^{n-2}}{q^n}+O_{Q,\nu,q}(t^{n-2-\kappa-\nu}).$$
\end{Cor}

\subsection{Second moment formula for congruence subgroups}
In view of the observation \eqref{equ:obs}, in order to apply the method developed in  \cite{KelmerYu2018b} to this setting, we need a formula for the second moment of the Siegel transform when the group $\SL_n(\Z)$ is replaced with an appropriate congruence subgroup.

Rogers' proof of his moment formula seems very special to the space of lattices, and it is not immediately clear how to generalize it to incorporate a congruence condition. 
In \cite{KelmerYu2018}, Kelmer and Yu introduced another method, relating the Siegel transform to incomplete Eisenstein series and relying on their spectral theory to obtain a second moment formula. This spectral approach is much more flexible and can be applied for many cases and in particular to congruence groups. The main technical result of this paper is thus an adaptation of this method to prove a second moment formula in this setting. We now give a special case of this formula that could be of independent interest.

Let $G=\SL_n(\R)$ and  $\G=\SL_n(\Z)$ and let $\mu$ denote the Haar measure of $G$ normalized so that $\mu(\G\bk G)=1$. Let $L\leq G$ denote the stabilizer of $\bm{e}_n=(0,\ldots,0,1)$ and note that there is a natural identification of $L\bk G$ with $\dot{\R}^n:=\R^n\setminus\{\bm{0}\}$ given by $Lg\mapsto \bm{e}_ng$.
 For $q\in \N$ we let $\G_1(q)\leq \G$ denote the congruence subgroup 
\begin{equation}
\G_1(q)=\{\g\in \G: \bm{e}_n\g\equiv \bm{e}_n\Mod{q}\}.
\end{equation}
 For a function on $f$ on $\R^n$ of sufficiently fast decay, we define its incomplete Eisenstein series by 
$$\Theta^{n,q}_f(g)=\sum_{\g\in \G_1(q)\cap L \bk \G_1(q)}f(\bm{e}_n\gamma g).$$
Note that $\Theta^{n,q}_f$ is a function on $\G_1(q)\bk G$. We then show the following formulas for the first and second moments (see \corref{cor:firstmoment} and \thmref{thm:secondmomentinein} for the general result).
 \begin{Thm}\label{thm:secondmomentIntro}
Let $n\geq 3$ be an integer. Let $f$ be a bounded compactly supported function on $L\bk G\cong \dot{\R}^n$.  Then
\begin{equation}\label{equ:firstmoment}
\int_{\G_1(q)\bk G}\Theta^{n,q}_{f}(g)d\mu(g)=\frac{\int_{\R^n}f(\bm{x})d\bm{x}}{\zeta(n)},
\end{equation}
and
\begin{equation}\label{equ:secondmomentIntro}
\int_{\G_1(q)\bk G}|\Theta^{n,q}_{f}(g)|^2d\mu(g)=\frac{|\int_{\R^n}f(\bm{x})d\bm{x}|^2}{q^n\zeta(n)\zeta_q(n)} +\frac{1}{\zeta(n)}\int_{\R^n}|f(\bm{x})|^2d\bm{x}+\frac{\delta_q}{\zeta(n)}\int_{\R^n}f(\bm{x})\overline{f(-\bm{x})}d\bm{x},
\end{equation}
where $\delta_q=1$ if $q\in \{1,2\}$ and $\delta_q=0$ if $q\geq 3$ and 
$\zeta_q(n)=\sum_{\substack{k\geq 1\\ \gcd(k,q)=1}}k^{-n}$.
\end{Thm}
 
\begin{rem}
Here the first moment formula \eqref{equ:firstmoment} also holds for $n=2$, see \corref{cor:firstmoment}. We also note that the orbit $\bm{e}_n\G_1(q)\subseteq \Z^n$ is precisely all primitive vectors in $\Z^n$ that are congruent to $\bm{e}_n$ modulo $q$, so the incomplete Eisenstein series  $\Theta_f^{n,q}$ is closely related to the Siegel transform when incorporating this congruence condition. If we want to change the congruence condition to $\bm u\equiv\bm p \Mod{q}$ for a different vector $\bm p\neq \bm{e}_n$, it can be done by replacing $\G_1(q)$ by an appropriate conjugate.
\end{rem}


\subsection{Irrational shifts}
The relation between values of an inhomogeneous form at integer points and values of the corresponding homogeneous form on integer points satisfying congruence conditions clearly no longer holds when the shift is irrational. Nevertheless, using the fact that the dependance on $q$ in our second moment formula is very explicit, we can still get some results for irrational shifts by approximating them by rational ones, as long as we have sufficient control on how fast the denominators grow. This control on the growth of the approximating vectors can be obtained for irrational vectors satisfying certain Diophantine conditions for which we can show the following.

\begin{Thm}\label{thm:irrational}
Let $n\geq 5$ and let $\bm{\alpha}\in \R^n$ be an irrational vector with Diophantine exponent $\omega_{\bm{\alpha}}$ and uniform Diophantine exponent $\hat{\omega}_{\bm{\alpha}}$ $($see \secref{sec:dio}$)$. Assume that $\omega_{\bm{\alpha}}<\infty$ and that $\hat{\omega}_{\bm{\alpha}}>\frac{2}{n-2}$. Then for any $\kappa\in \left(0,\frac{(n-2)\hat{\omega}_{\bm{\alpha}}-2}{n\left(1+\omega_{\bm{\alpha}}+\omega_{\bm{\alpha}}\hat{\omega}_{\bm{\alpha}}\right)}\right)$, for any $\xi\in \R$, and for almost every non-degenerate indefinite quadratic form $Q$ in $n$ variables, the system of inequalities \eqref{equ:maininequ} has integer solutions for all sufficiently large $t$.
\end{Thm}
\begin{rem}
We note that here the assumption that $n\geq 5$ is to ensure the existence of $\bm{\alpha}$ satisfying the above two Diophantine conditions (see \rmkref{rmk:setdi}), and in fact we can handle a slightly larger set of irrational vectors (see \thmref{thm:moregeneral}). We also note that our result for irrational shift is much weaker (for example our exponent is always smaller than $\hat{\omega}_{\bm{\alpha}}$ which is less or equal to $1$ (see \rmkref{rmk:setdi}) while for generic shifts the correct critical exponent is $n-2$). This could be just an artefact of our proof, but it might suggest that there are special irrational shifts that behave much worse than the generic ones, and it reflects the sensitivity of the Diophantine problem under consideration.     
\end{rem}
\subsection{Further generalizations}
We have not addressed here the uniform results for approximating all targets simultaneously 
as done in \cite{Bourgain2016,GhoshKelmer2018,KelmerYu2018b}
for homogeneous forms. We note that while there are some obstacles getting a uniform analogue of \thmref{thm:irrational}, our method can be easily adapted to deal with this type of problem in the fixed rational shift setting. We refer the reader to the proof of Corollary 4 of \cite{KelmerYu2018b} to see how such a result follows from the second moment formula.

We also note that, just as in \cite{KelmerYu2018b}, one can extend these methods to give similar results not only to quadratic forms but to shifts (either rational or irrational and satisfying similar Diophantine conditions) of more general homogeneous polynomials of higher degrees.

%
%
%
\subsection*{Notation}
Let us fix some notation throughout this paper. For two positive quantities $A$ and $B$, we will use the notation $A\ll B$  {or $A=O(B)$} to mean that there is a constant $c>0$ such that $A\leq cB$, and we will use subscripts to indicate the dependence of the constant on parameters. We will write $A\asymp B$ for $A\ll B\ll A$. For any $t>0$ we denote by $B_t\subset \R^n$ the open Euclidean ball centered at the origin with radius $t$. For any measurable function $f$ on $\R^n$ we denote by $\vol(f):=\int_{\R^n}f(\bm{x})d\bm{x}$. 
\subsection*{Acknowledgements}
The authors would like to thank Jens Marklof for his comments and for bringing their attention to some references. The authors would also like to thank an anonymous referee for a thoughtful report, especially for suggesting an alternative proof to the second moment formula in \thmref{thm:moment}. 
\section{Moment formulas of incomplete Eisenstein series
}\label{sec:incomplete}
In this section we prove the main technical result of this paper, which is a more general second moment formula than \eqref{equ:secondmomentIntro} for certain translates of the incomplete Eisenstein series $\Theta_f^{n,q}$ defined in the introduction. To prove such a moment formula we need to work on homogeneous spaces of different ranks, and we thus keep this parameter in our notation. Let $n\geq 2$. Let $G_n=\SL_n(\R)$ and let $\G_n=\SL_n(\Z)$. Let $\mu_n$ be a Haar measure of $G_n$ normalized such that $\mu_n(\G_n\bk G_n)=1$. 

Let $P_n\leq G_n$ be the identity component of the maximal parabolic subgroup of $G_n$ fixing the line spanned by $\bm{e}_n=(0,\ldots,0,1)\in\R^n$ under the right multiplication action, and let $L_n\leq P_n$ be the stabilizer of $\bm{e}_n$. Since this action on $\dot{\R}^n=\R^n\setminus\{\bm{0}\}$ is transitive, it induces an identification between the homogeneous space $L_n\bk G_n$ and $\dot{\R}^n$ identifying $L_ng$ with $\bm{e}_ng$, the bottom row of $g$.

Let $\G\leq \G_n$ be a finite-index subgroup of $\G_n$. Given a bounded compactly supported function $f$ on $L_n\bk G_n$ (that we think of as a left $L_n$-invariant function on $G_n$),
the corresponding \textit{incomplete Eisenstein series at $P_n$}, denoted by $\Theta^{\G}_f$ on $\G\bk G_n$ is defined by
$$\Theta^{\G}_f(g)=\sum_{\gamma\in \G\cap P_n\bk \G}f(\gamma g).$$
We note that it is not difficult to check that $\G_n\cap P_n=\G_n\cap L_n$. Since $\G\leq \G_n$ we also have $\G\cap P_n=\G\cap L_n$. Together with the assumption that $f$ is left $L_n$-invariant, this implies that the series defining $\Theta^{\G}_f$ is well-defined. Moreover, it is easy to check that $\Theta_f^{\G}$ is left $\G$-invariant, and since $f$ is bounded and compactly supported, the series for $\Theta_f^{\G}(g)$ is actually a finite sum for any $g\in G_n$, and hence absolutely converges.

For any integer $q\in \N$, let
$$\G_1^n(q):=\left\{\gamma\in \G_n\ |\ \bm{e}_n\gamma\equiv \bm{e}_n \Mod{q}\right\}$$  
be the congruence subgroup consisting of elements in $\G_n$ whose bottom row is congruent to $\bm{e}_n$ modulo $q$. We denote by $X_{n,q}=\G_1^n(q)\bk G_n$  the corresponding homogeneous space. For any $\ell\in \Z$ we denote by $[\ell]$ its reduction in $\Z/q\Z$, and denote by $[\bar{\ell}]$ the inverse of $[\ell]$ in $(\Z/q\Z)^{\times}$ if $\gcd (\ell,q)=1$. For each $[\ell]\in (\Z/q\Z)^{\times}$ since the action of $\G_n$ on the set of primitive integer vectors, $\Z^n_{\rm pr}$, is transitive, we can take $\tau_{[\ell]}\in \G_n$ such that 
\begin{equation}\label{equ:tau}
\bm{e}_n\tau_{[\ell]}\equiv [\bar{\ell}]\bm{e}_n \Mod{q}.
\end{equation}
We note that the choice of $\tau_{[\ell]}$ is not unique and \eqref{equ:tau} implies that $\tau_{[\ell]}$ normalizes $\G_1^n(q)$. For simplicity of notation, when $\G=\G_1^n(q)$ we abbreviate $\Theta_f^{\G_1^n(q)}$ by $\Theta_f^{n,q}$. Since $\G_1^n(q)\cap P_n=\G_1^n(q)\cap L_n$, this definition of $\Theta_f^{n,q}$ here coincides with the definition given in the introduction. We prove the following inner product formula for $\Theta_f^{n,q}$.
\begin{Thm}\label{thm:secondmomentinein}
Let $f_1$ and $f_2$ be any two bounded and compactly supported functions on $L_n\bk G_n\cong \dot{\R}^n$. Then for any $q\geq 1$, $n\geq 2$ and $[\ell]\in (\Z/q\Z)^{\times}$ we have
\begin{align}\label{equ:secondmoment}
&\int_{X_{n,q}}\Theta^{n,q}_{f_1}(g)\Theta_{f_2}^{n,q}(\tau_{[\ell]}g)d\mu_n(g)-\frac{1}{\zeta(n)}\left(\delta_{[\ell][1]}\vol(f_1f_2)+\delta_{[\ell][-1]}\vol(f_1\tilde{f}_2)\right)\\
&=\left\{\begin{array}{ll}
\frac{\vol(f_1)\vol(f_2)}{q^n\zeta(n)\zeta_q(n)}, & n\geq 3\\
\sum_{w\neq 0}\frac{\varphi(w_q)}{qw_q\zeta(2)}\int_{\R^3}f_1\left(x_1,x_2\right)f_2\left(\frac{qwx_2}{x_1^2+x_2^2}+tx_1,\frac{-qwx_1}{x_1^2+x_2^2}+tx_2\right)dtdx_1dx_2, & n=2, \nonumber
\end{array}\right.
\end{align}
where $\tilde{f}(\bm{x}):=f(-\bm{x})$, $\delta_{[\ell][1]}$ is the Kronecker delta function, $\zeta_q(n)=\sum_{\substack{k\geq 1\\ \gcd(k,q)=1}}k^{-n}$, $\varphi$ is the Euler's totient function and for any integer $w\neq 0$, $w_q$ is the largest positive divisor of $w$ such that $\gcd(w_q,q)=1$.
\end{Thm}
\begin{rmk}
When $n=2$ the inner product formula in \eqref{equ:secondmoment} generalizes a recent result in \cite[Theorem 2.1]{KleinbockYu2020}. This formula will not be needed for our application, but it is of independent interest.
\end{rmk}
The rest of this section is devoted to proving the inner product formula \eqref{equ:secondmoment}. In \cite{KelmerYu2018}, a similar formula was proved for the symplectic group, by relating the incomplete Eisenstein series as an integral over Eisenstein series and using the analytic properties of the constant term of the Eisenstein series. 
While the same approach can also be applied here (and was done in an earlier version of this paper), we will use instead a more direct approach using unfolding and some elementary properties of the space of lattices.

Before proceeding to the proof, let us first make a few remarks about the integral in the left hand side of \eqref{equ:secondmoment}. First we note that since $\tau_{[\ell]}$ normalizes $\G_1^n(q)$, the function $\G_1^n(q)g\mapsto \Theta_f^{n,q}(\tau_{[\ell]}g)$ is well-defined on $X_{n,q}$ and so is the above integral. Moreover, suppose $\tau_{[\ell]}'\in \G_n$ is another element in $\G_n$ satisfying \eqref{equ:tau}, then by definition we have $\tau_{[\ell]}'\tau_{[\ell]}^{-1}\in \G_1^n(q)$. Using the left $\G_1^n(q)$-invariance of $\Theta_{f}^{n,q}$, we have 
$$\Theta_f^{n,q}(\tau_{[\ell]}g)=\Theta_{f}^{n,q}(\tau_{[\ell]}'g)$$ 
for any $g\in G_n$, implying that the above integral is independent of the choice of $\tau_{[\ell]}$. Moreover, for later reference we note that for any $[\ell_1],[\ell_2]\in(\Z/q\Z)^{\times}$ we have
\begin{equation}\label{equ:choiceoftau}
\Theta_f^{n,q}(\tau_{[\ell_1]}\tau_{[\ell_2]}g)=\Theta_{f}^{n,q}(\tau_{[\ell_1\ell_2]}g).
\end{equation}
Finally we note that the condition that $f$ being bounded and compactly supported can be relaxed to the condition that $f$ is bounded and nonnegative, see \cite[Remark 5.13]{KelmerYu2018}. In particular, we can apply \eqref{equ:secondmoment} to indicator functions of any finite-volume subsets of $\R^n$.
\subsection{\textbf{Coordinates and measures}}
Let $P_n$ and $L_n$ be as above. There is a Langlands decomposition $P_n=U_nA_nM_n$ with 
$$U_n=\left\{u_{\bm{t}}=\begin{pmatrix}
I_{n-1} & \bm{t}^t\\
\bm{0} & 1\end{pmatrix}\ \left|\ \bm{t}=(t_1,\ldots,t_{n-1})\in \R^{n-1}\ 
\right.\right\}
,$$
$$A_n= \left\{a_{y}=\textrm{diag}(y^{\frac{1}{n-1}},\ldots, y^{\frac{1}{n-1}}, y^{-1})\ \left|\ y>0 \right.\right\}
$$
and 
$$
M_n=\left\{ \widetilde{m}=\begin{pmatrix}
m& \bm{0}^t \\
\bm{0} & 1
\end{pmatrix}\ \left|\  m\in G_{n-1}\right.\right\}.$$

We note that $L_n=U_nM_n$. Fix a maximal compact subgroup $K_n=\SO_n(\R)$ and with slight abuse of notation, we denote by $K_{n-1} :=M_n\cap K_n$. Note that $K_{n-1}\cong \SO_{n-1}(\R)$ is a maximal compact subgroup of $M_n\cong G_{n-1}$, and we can identify $K_{n-1}\bk K_n$ with the unit sphere $S^{n-1}:=\{\bm{x}\in \R^n\ |\ \|\bm{x}\|=1\}$ by identifying $K_{n-1}k$ with $\bm{e}_nk$.

Recall the identification between $L_n\bk G_n$ and $\dot{\R}^n$ sending $L_ng$ to $\bm{e}_ng$. Let $d\bm{x}$ be the Lebesgue measure on $\R^n$ that we think of as a measure on $L_n\bk G_n$. We also use the following polar coordinates on $\dot{\R}^n\cong L_n\bk G_n$: as a consequence of the Iwasawa decomposition for $G_n$ we have $G_n=U_nM_nA_nK_n$. Thus any element in $L_n\bk G_n$ can be represented uniquely by $a_yk$ with $a_y\in A_n$ and $k\in K_{n-1}\bk K_n\cong S^{n-1}$. In these coordinates we have that 
\begin{equation}\label{equ:polar}
d\bm{x}(a_yk)=\frac{2\pi^{n/2}}{\G(\frac{n}{2})}\frac{dy}{y^{n+1}}d\sigma_n(k),
\end{equation}
where $d\sigma_n(k)$ is the probability right $K_n$-invariant measure on $K_{n-1}\bk K_n$ (which is the surface measure on $S^{n-1}$) and $\G(s)=\int_0^{\infty}t^{s-1}e^{-t}dt$ for $\Re(s)>0$ is the Gamma function.

Next, let $L_n(\Z)=L_n\cap \G_n$ and denote by $\mu_{L_n}$ the probability Haar measure on $L_n(\Z)\bk L_n$. Explicitly, for any $h\in L_n$ writing $h=u_{\bm{t}}\tilde{m}$ with $u_{\bm{t}}\in U_n$, $\tilde{m}\in M_n$ and identifying $M_n$ with $G_{n-1},$ we have
\begin{equation}\label{equ:measurel}
d\mu_{L_n}(h)=d\mu_{L_n}(u_{\bm{t}}\tilde{m})=d\bm{t}d\mu_{n-1}(m),
\end{equation}
where when $n=2$, $\mu_{n-1}$ is the probability measure on the trivial group. Using this decomposition, for any smooth and compactly supported function $F$ on $G_n$ we have
\begin{equation}\label{equ:haar}
\int_{G_n}F(g)d\mu_n(g)=\omega_n\int_{L_n\bk G_n}\int_{L_n}F(ha_yk)d\mu_{L_n}(h)\frac{dy}{y^{n+1}}d\sigma_n(k),
\end{equation}
where $\omega_n=\frac{2}{\xi(n)}$ with $\xi(s)=\pi^{-s/2}\G(\frac{s}{2})\zeta(s)$ the Riemann Xi function.
\subsection{\textbf{More on the congruence subgroup $\G_1^n(q)$}}
Let $\nu_{n,q}=\mu_n(X_{n,q})=[\G_n: \G_1^n(q)]$ be the index of $\G_1^n(q)$ in $\G_n$.  
We first prove a formula for $\nu_{n,q}$. 
\begin{Lem}\label{lem:index}
For any $n\geq 2$ and $q\in \N$ we have
$$\nu_{n,q}=\frac{q^n\zeta_q(n)}{\zeta(n)}.$$
\end{Lem}
\begin{proof}
Define $\G_0^n(q)< \G_n$ such that 
$$\G_0^n(q)=\left\{\gamma\in \G_n\ |\ \bm{e}_n\gamma\equiv \ell\bm{e}_n\Mod{q}\ \textrm{for some $[\ell]\in (\Z/q\Z)^{\times}$}\right\}.$$
By  \cite[Proposition 2.1]{BazHuangLee2019} we have $[\G_n : \G_0^n(q)]=q^{n-1}\prod_{p| q}\frac{1-p^{-n}}{1-p^{-1}}$. Consider the map
$$f: \G_0^n(q)\to (\Z/q\Z)^{\times}$$
sending $\gamma=(a_{ij})_{n\times n}\in \G_0^n(q)$ to $[a_{nn}]\in(\Z/q\Z)^{\times}$. It is easy to check that $f$ is a surjective group homomorphism with the kernel given by $\G_1^n(q)$, so that 
$$[\G_0^n(q): \G_1^n(q)]=\left|( \Z/q\Z)^{\times}\right|=q\prod_{p| q}(1-p^{-1}).$$
We thus have
$$\nu_{n,q}=[\G_n:\G_1^n(q)]=[\G_n:\G_0^n(q)][\G_0^n(q): \G_1^n(q)]=
q^n\prod_{p| q}(1-p^{-n})=\frac{q^n\zeta_q(n)}{\zeta(n)}$$
concluding the proof. 
\end{proof}
For any $(\bm{p},q)\in \Z^n\times \N$ with $\gcd (\bm{p}, q)= 1$, we define
$$\cO(\bm{p}):=\left\{\bm{v}\in \Z^n\setminus \{\bm{0}\}\ |\ \bm{v}\equiv \bm{p} \Mod{q}\right\},$$
and for any $[\ell]\in (\Z/q\Z)^{\times}$ we define
$$\cO(\bm{p};[\ell]):=\left\{\bm{v}\in \Z^n_{\rm pr}\ |\ \bm{v}\equiv \bar{[\ell]}\bm{p}\Mod{q}\right\}.$$
The following lemma gives an identification between the coset representatives for $\G_1^n(q)\cap P_n\bk \G_1^n(q)$ and $\cO(\bm{e}_n,[1])$, the set of primitive integer vectors congruent to $\bm{e}_n$ modulo $q$.
\begin{Lem}\label{lem:id2}
The map from $\G_1^n(q)$ to $\Z_{\rm pr}^n$ sending $\gamma \in \G_1^n(q)$ to $\bm{e}_n\gamma$ induces an identification between $\G_1^n(q)\cap P_n\bk \G_1^n(q)$ and $\cO(\bm{e}_n; [1])$. More generally, the map sending  $\gamma\tau_{[\ell]}$ to $\bm{e}_n\gamma\tau_{[\ell]}$ gives a bijection between
 $\G_1^n(q)\cap P_n\bk \G_1^n(q)\tau_{[\ell]}$ and  $\cO(\bm{e}_n; [\ell])$.
 \end{Lem}
\begin{proof}
For the first identification, it suffices to show that $\bm{e}_n\G_1^n(q)=\cO(\bm{e}_n,[1])$ and that the stabilizer of $\bm{e}_n$ in $\G_1^n(q)$ is $\G_1^n(q)\cap P_n$. For the first claim, the containment $\bm{e}_n\G_1^n(q)\subset \cO(\bm{e}_n,[1])$ is clear from the definition. For the other containment, since the right multiplication action of $\G_n$ on $\Z^n_{\rm pr}$ is transitive, for any $\bm{v}\in \cO(\bm{e}_n; [1])$, there exists some $\gamma \in \G_n$ such that $\bm{e}_n\gamma =\bm{v}$. Hence $\bm{e}_n\gamma= \bm{v}\equiv \bm{e}_n\Mod{q}$ implying that $\gamma\in \G_1^n(q)$. For the second claim, since $L_n$ is the stabilizer of $\bm{e}_n$, we have the stabilizer of $\bm{e}_n$ in $\G_1^n(q)$ is $\G_1^n(q)\cap L_n=\G_1^n(q)\cap P_n$, proving the claim. The second part follows from the first identification and the relation $\cO(\bm{e}_n; [\ell])=\cO(\bm{e}_n; [1])\tau_{[\ell]}$.
\end{proof}

\subsection{Unfolding}
In this subsection, we prove some preliminary identities using the standard unfolding trick, captured in the following. We note that $\G_1^n(q)$ satisfies the below assumptions on the lattice $\G$.
\begin{Lem}\label{lem:unfolding}
Let $\G<\G_n$ be a finite-index subgroup of $\G_n$ satisfying that $\G\cap P_n=\G_n\cap P_n$. For any bounded function $f$ on $L_n\bk G_n$ satisfying that $|f(a_yk)|\ll y^\sigma$ for some $\sigma>n$, and $\Psi\in L^2(\G\bk G_n)$ 
we have 
\begin{equation}
\int_{\G\bk G_n}
\Theta_f^\G(g)\Psi(g)d\mu_n(g)=\omega_n \int_{L_n\bk G_n}f(a_y k)\cP(\Psi)(a_yk)\frac{dy}{y^{n+1}}d\sigma_n(k),
\end{equation}
where $\cP$ is a period operator sending a function $\Psi$ on $\G\bk G_n$ to a function on $L_n\bk G_n$ given by 
 \begin{equation}
\cP(\Psi)(a_y k):=\int_{L_n(\Z)\bk L_n}\Psi(ha_y k)d\mu_{L_n}(h).
\end{equation}
\end{Lem}
\begin{proof}
First note that since $\G_n\cap P_n=\G_n\cap L_n=L_n(\Z)$ we also have that $\G\cap P_n=L_n(\Z)$. Now the condition $f(a_yk)\ll y^\sigma$ with $\sigma>n$ implies that the series 
$$\Theta^\G_f(g)=\sum_{\g\in L_n(\Z)\bk \G}f(\g g)$$
absolutely converges and hence the standard unfolding argument gives
\begin{align*}
\int_{\G\bk G_n}\Theta_f^\G(g)\Psi(g)d\mu_n(g)&=\int_{L_n(\Z)\bk G_n}f(g)\Psi(g)d\mu_n(g)\\
&=\omega_n\int_{L_n\bk G_n}\int_{L_n(\Z)\bk L_n} f(a_yk)\Psi(ha_y k)d\mu_{L_n}(h)\frac{dy}{y^{n+1}}d\sigma_n(k)\\
&= \omega_n\int_{L_n\bk G_n}f(a_yk)\cP(\Psi)(a_y k)\frac{dy}{y^{n+1}}d\sigma_n(k),
\end{align*}
where the second line follows from the decomposition \eqref{equ:haar} of the Haar measure $\mu_n$ and the left $L_n$-invariance of $f$.
\end{proof}

Using this with $\Psi=1$ gives the following first moment formula.
\begin{Cor}\label{cor:firstmoment}
Let $\G<\G_n$ be as in \lemref{lem:unfolding}. For a function $f$ on $L_n\bk G_n\cong  \dot \R^n$ as in \lemref{lem:unfolding} we have
\begin{equation}\label{equ:firstincomp}
\int_{\G\bk G_n}\Theta^{\G}_{f}(g)d\mu_n(g)=\omega_n\int_{L_n\bk G_n}f(a_yk)\frac{dy}{y^{n+1}}d\sigma_n(k)=\frac{\int_{\R^n}f(\bm{x})d\bm{x}}{\zeta(n)}.
\end{equation}
\end{Cor}
\subsection{Proof of \thmref{thm:secondmomentinein}}
In this subsection we give the proof of \thmref{thm:secondmomentinein}. Let $f$ be a bounded and compactly supported function. For any $[\ell]\in (\Z/q\Z)^{\times}$ denote by $\Theta_f^{n,q,\tau_{[\ell]}}$ the incomplete Eisenstein series on $\G_1^n(q)$ twisted by $\tau_{[\ell]}$, that is, for any $g\in G_n$, $\Theta_f^{n,q,\tau_{[\ell]}}(g):=\Theta_f^{n,q}(\tau_{[\ell]}g)$.  In view of \lemref{lem:unfolding}, to prove \thmref{thm:secondmomentinein} it suffices to compute the period function $\cP\left(\Theta_f^{n,q,\tau_{[\ell]}}\right)$. A similar computation was done in \cite{KelmerYu2018} for symplectic groups by using some general theory of Eisenstein series. However, as we mentioned above, in this situation there is a more direct way of computing this period function (following the argument in  \cite[Lemma 7.7]{MarklofStrombergsson2010}). We prove the following formula for $\cP\left(\Theta_f^{n,q,\tau_{[\ell]}}\right)$ which, together with  \lemref{lem:unfolding} directly implies \thmref{thm:secondmomentinein}.
\begin{Prop}\label{prop:elementary}
Let $f :\R\to \C$ be a bounded and compactly supported function. For any $q\geq 1, n\geq 2$ and $[\ell]\in (\Z/n\Z)^{\times}$ let $\cP\left(\Theta_f^{n,q,\tau_{[\ell]}}\right)$ be as in \eqref{equ:Pf}. Then 
for $n\geq 3$ 
\begin{eqnarray*}
\cP\left(\Theta_f^{n,q,\tau_{[\ell]}}\right)(\bm{x})&=&\delta_{[\ell][1]}f(\bm{x})+\delta_{[\ell][-1]}f(-\bm{x})+\frac{\vol(f)}{q^n\zeta_q(n)},
\end{eqnarray*}
while for $n=2$
\begin{eqnarray*}\cP\left(\Theta_f^{n,q,\tau_{[\ell]}}\right)(\bm{x})&=&\delta_{[\ell][1]}f(\bm{x})+\delta_{[\ell][-1]}f(-\bm{x})\\
&&+\sum_{w\neq 0}\frac{\varphi(w_q)}{qw_q}\int_{\R}f\left(\frac{qwx_2}{x_1^2+x_2^2}+tx_1,\frac{-qwx_1}{x_1^2+x_2^2}+tx_2\right)dt.
\end{eqnarray*}
where $\zeta_q(n)$, $\varphi$ and $w_q$ are as in \thmref{thm:secondmomentinein}.
\end{Prop}
\begin{proof}
We first give a more explicit description of the period operator $\cP$. Let $M_n(\Z)=M_n\cap \G_n$ and $U_n(\Z)=U_n\cap \G_n$. Using the measure decomposition \eqref{equ:measurel} and the relation $U_n(\Z)M_n(\Z)=L_n(\Z)$ and identifying $L_n\bk G_n$ with $\dot{\R}^n$ we have
\begin{equation}\label{equ:Pf}
\cP\left(\Theta_f^{n,q,\tau_{[\ell]}}\right)(\bm{e}_na_yk)=\int_{M_n(\Z)\bk M_n}\int_{U_n(\Z)\bk U_n}\Theta^{n,q}_f(\tau_{[\ell]}u_{\bm{t}}\widetilde{m}a_yk)d\bm{t}d\mu_{n-1}(m).
\end{equation}

Let $\mu: \N\to \{0,\pm 1\}$ be the M$\ddot{\textrm{o}}$bius function. For any $\bm{v}=(v_1,\ldots, v_n)\in \R^n$, we denote by $\bm{v}^{M_n}=(v_1,\ldots, v_{n-1})\in \R^{n-1}$. Firstly, by the second part of \lemref{lem:id2}, we can rewrite
\begin{align*}
\Theta_f^{n,q}(\tau_{[\ell]}g)&=\sum_{\bm{v}\in \cO(\bm{e}_n; [\ell])}f(\bm{v}g)\\
&=\delta_{[\ell][1]}f(\bm{e}_ng)+\delta_{[\ell][-1]}f(-\bm{e}_ng)+\mathop{\sum_{\bm{v}\equiv [\bar{\ell}]\bm{e}_n (\mod q)}}_{\bm{v}^{M_n}\neq \bm{0}}\sum_{d| \gcd \left(\bm{v}\right)}\mu(d)f(\bm{v}g)\\
&=\delta_{[\ell][1]}f(\bm{e}_ng)+\delta_{[\ell][-1]}f(-\bm{e}_ng)+\mathop{\sum_{d\geq 1}}_{\gcd (d,q)=1}\mu(d)\mathop{\sum_{\bm{v}\equiv [\bar{d}\bar{\ell}]\bm{e}_n (\mod q)}}_{\bm{v}^{M_n}\neq \bm{0}}f(d\bm{v}g),
\end{align*}
where for the second equality we used the fact that $\pm \bm{e}_n\in \cO(\bm{e}_n; [\ell])$ if and only if $[\ell]=[\pm1]$ and the identity $\sum_{d| n}\mu(d)=\delta_{n1}$, and for the third equality we used that $\gcd \left(\bm{v}\right)$ is coprime to $q$ and did a change of variable $\bm{v}\mapsto d\bm{v}$. For any $\lambda>0$ denote by $f_{\lambda}(\bm{x}):=f(\lambda \bm{x})$. Taking $g=u_{\bm{t}}\widetilde{m}a_yk$ and identifying $U_n(\Z)\bk U_n$ with $[0,1)^{n-1}$ we get  that 
\begin{eqnarray*}
\mathcal P^{n,q,\tau_{[\ell]}}_f(\bm{e}_na_yk)&=& \left(\delta_{[\ell][1]}f(\bm{e}_na_yk)+\delta_{[\ell][-1]}f(-\bm{e}_na_yk)\right)\\
&&+
\mathop{\sum_{d\geq 1}}_{\gcd (d,q)=1}\mu(d)\int_{M_n(\Z)\bk M_n}\int_{[0,1)^{n-1}}\mathop{\sum_{\bm{v}\equiv [\bar{d}\bar{\ell}]\bm{e}_n (\mod q)}}_{\bm{v}^{M_n}\neq \bm{0}}f_d(\bm{v}u_{\bm{t}}\widetilde{m}a_yk)d\bm{t}d\mu_{n-1}(m).
\end{eqnarray*}
It thus suffices to compute the inner integral
\begin{equation}\label{e:cI}
\cI:=\int_{M_n(\Z)\bk M_n}\int_{[0,1)^{n-1}}\mathop{\sum_{\bm{v}\equiv [\bar{d}\bar{\ell}]\bm{e}_n (\mod q)}}_{\bm{v}^{M_n}\neq \bm{0}}f_d(\bm{v}u_{\bm{t}}\widetilde{m}a_yk)d\bm{t}d\mu_{n-1}(m).
\end{equation}
Fix $\bar{d}\bar{\ell}\in \Z$ a lift of $[\bar{d}\bar{\ell}]$ in $(\Z/ q\Z)^{\times}$. Note that any $\bm{v}\in \Z^n$ satisfying $\bm{v}\equiv [\bar{d}\bar{\ell}]\bm{e}_n (\mod q)$ and $\bm{v}^{M_n}\neq \bm{0}$ can be written uniquely as $\bm{v}=q\bm{w}+\bar{d}\bar{\ell}\bm{e}_n$ for some $\bm{w}\in \Z^n$ with $\bm{w}^{M_n}\neq \bm{0}$. Moreover, in this case by direct computation we have
$$ \bm{v}u_{\bm{t}}\widetilde{m}=\left(q\bm{w}+\bar{d}\bar{\ell}\bm{e}_n\right)u_{\bm{t}}\widetilde{m}= \left(q\bm{w}^{M_n}m, qw_n+ \bar{d}\bar{\ell}+ q(w_1t_1+\ldots+w_{n-1}t_{n-1})\right),$$
so that the set 
$$\left\{\bm{v}u_{\bm{t}}\widetilde{m}\ |\  \bm{v}\equiv [\bar{d}\bar{\ell}]\bm{e}_n (\mod q) ,\; \bm{v}^{M_n}\neq \bm{0}\right\}$$
can be written explicitly as 
$$\{\left(q\bm{w}^{M_n}m, qw_n+ \bar{d}\bar{\ell}+ q(w_1t_1+\ldots+w_{n-1}t_{n-1})\right)\ |\  \bm{w}\in \Z^n,\; \bm{w}^{M_n}\neq \bm{0}\}.$$

Since $\bm{w}^{M_n}\neq \bm{0}$, there exists some $w_i\neq 0$ for some $1\leq i\leq n-1$.
For this fixed $i$, note that as $w_n$ runs through all the integers, the intervals 
$$[qw_n+\bar{d}\bar{\ell}+q\sum_{j\neq i}w_jt_j, qw_n+\bar{d}\bar{\ell}+q\sum_{j\neq i}w_jt_j+qw_i)$$ 
cover $\R$ exactly $|w_i|$ times, implying that for a fixed $\bm{w}^{M_n}\neq \bm{0}$ we have
\begin{align*}
&\int_{[0,1)^{n-1}}\sum_{w_n\in \Z}f_d\left(\left(q\bm{w}^{M_n}m, qw_n+ \bar{d}\bar{\ell}+ q(w_1t_1+\ldots+w_{n-1}t_{n-1})\right)a_yk\right)d\bm{t}\\
&=\int_{[0,1)^{n-1}}|w_i|\int_{\R}f_d\left((q\bm{w}^{M_n}m, qw_i t_i)a_yk\right)dt_i\prod_{\substack{1\leq j\leq n-1\\ j\neq i}}dt_j=\int_{\R}f_{qd}((\bm{w}^{M_n}m, t)a_yk)dt,\\
\end{align*}
where for the second equality we did a change of variable $w_it_i\mapsto t$.
Summing over all $\bm{w}^M_n\in \Z^{n-1}\setminus\{\bm{0}\}$ we get that  
\begin{eqnarray*}
\int_{[0,1)^{n-1}}\mathop{\sum_{\bm{v}\equiv [\bar{d}\bar{\ell}]\bm{e}_n (\mod q)}}_{\bm{v}^{M_n}\neq \bm{0}}f_d(\bm{v}u_{\bm{t}}\widetilde{m}a_yk)d\bm{t}&=&\sum_{\bm{w}^{M_n}\neq \bm{0}} \int_{\R}f_{qd}((\bm{w}^{M_n}m, t)a_yk)dt\\
&=&\sum_{\bm{w}^{M_n}\neq \bm{0}}F_{f_{qd}}\left(\bm{w}^{M_n}m\right),
\end{eqnarray*}
where, for fixed $a_yk\in G_n$, we let $F_f: \R^{n-1}\to \C$ be defined by $F_f(\bm{z}):=\int_{\R}f\left((\bm{z},t)a_yk\right)dt$ (note that since $f$ is bounded and compactly supported, so is $F_f$). Plugging this back into \eqref{e:cI} we get that
\begin{align*}
\cI&=\int_{M_n(\Z)\bk M_n}\sum_{\bm{w}^{M_n}\neq \bm{0}}F_{f_{qd}}\left(\bm{w}^{M_n}m\right)d\mu_{n-1}(m).
\end{align*}
Now we separate the proof into two cases. First we assume that $n\geq 3$, then using the isomorphism $\G_{n-1}\bk G_{n-1}\cong M_n(\Z)\bk M_n(\Z)$, Siegel's Mean Value theorem \cite{Siegel1945} (noting that $\sum_{\bm{w}^{M_n}\neq \bm{0}}F_{f_{qd}}(\bm{w}^{M_n}m)$ is a Siegel transform defined on $M_n(\Z)\bk M_n$, see \secref{sec:Siegel} for more details on Siegel transforms) and the fact that $a_yk$ is of determinant one, we get
$$\cI=\int_{\R^{n-1}}F_{f_{qd}}\left(\bm{z}\right)d\bm{z}=\int_{\R^{n-1}}\int_{\R}f_{qd}\left((\bm{z},t)a_yk\right)d\bm{z}dt=\frac{\vol(f)}{q^nd^n}.$$
Thus using the identity $\sum_{\substack{d\geq 1 \\ \gcd (d,q)=1}}\frac{\mu(d)}{d^n}=\zeta_q(n)^{-1}$, we get for $\bm{x}=\bm{e}_na_yk$ 
$$\mathcal P^{n,q,\tau_{[\ell]}}_f(\bm{x})-\left(\delta_{[\ell][1]}f(\bm{x})+\delta_{[\ell][-1]}f(-\bm{x})\right)=\sum_{\substack{d\geq 1\\ \gcd(d,q)=1}}\mu(d)\frac{\vol(f)}{q^nd^n}=\frac{\vol(f)}{q^n\zeta_q(n)},
$$ as claimed. 

Next, when $n=2$, $M_n(\Z)\bk M_n$ is just a single point and $\mu_{n-1}$ is the probability measure supported on this point. Hence we get
$$\cI=\sum_{w\neq 0}F_{f_{qd}}(w)=\sum_{w\neq 0}\int_{\R}f_{qd}((w,t)a_yk)dt=\sum_{\substack{w\neq 0\\ d | w}}\frac{1}{qd}\int_{\R}f\left((qw, t)a_yk\right)dt,$$
where for the last equality we did change of variables $qdt\mapsto t$ and $dw\mapsto w$. Thus we have for $\bm{x}=\bm{e}_2a_yk$, $\mathcal P^{n,q,\tau_{[\ell]}}_f(\bm{x})-\left(\delta_{[\ell][1]}f(\bm{x})+\delta_{[\ell][-1]}f(-\bm{x})\right)$ equals
\begin{align*}
&\sum_{\substack{d\geq 1\\ \gcd(d,q)=1}}\mu(d)\sum_{\substack{w\neq 0\\ d| w}}\frac{1}{qd}\int_{\R}f\left((qw, t)a_yk\right)dt=\sum_{w\neq 0}\sum_{\substack{d\geq 1\\ d|w_q}}\frac{\mu(d)}{qd}\int_{\R}f\left((qw, t)a_yk\right)dt\\
&=\sum_{w\neq 0}\frac{\varphi(w_q)}{qw_q}\int_{\R}f\left((qw, t)a_yk\right)dt,
\end{align*}
where for the first equality we used the fact that for $d\geq 1$ satisfying $\gcd(d,q)=1$, the condition $d| w$ is equivalent to $d|w_q$, and for the second equality we used the identity $\sum_{\substack{d\geq 1\\ d| w_q}}\mu(d)/d=\varphi(w_q)/w_q$.
Finally, we conclude the proof by noting that for $(x_1,x_2)=(0,1)a_yk_{\theta}$ with $a_y=\left(\begin{smallmatrix}
y & 0\\
0 & y^{-1}\end{smallmatrix}\right)$ and $k_{\theta}=\left(\begin{smallmatrix}
\cos\theta & \sin\theta\\
-\sin\theta & \cos\theta\end{smallmatrix}\right)$, 
\begin{displaymath}
(qw,t)a_yk_{\theta}=\left((qw,0)+(0,t)\right)a_yk_{\theta}=\left(\frac{qwx_2}{x_1^2+x_2^2}+tx_1,\frac{-qwx_1}{x_1^2+x_2^2}+tx_2\right).\qedhere
\end{displaymath}
\end{proof}

\section{Moment formulas of Siegel transforms}\label{sec:Siegel}
Let $n\geq 2$ and let $AX_n=\rm{ASL}_n(\Z)\bk \rm{ASL}_n(\R)$ be the space of affine unimodular lattices in $\R^n$. For any bounded and compactly supported function $f: \R^n\to \C$ recall the \textit{Siegel transform}, $\hat{f} : AX_n \to \C$, is defined such that for any $\Lambda\in AX_n$
$$\hat{f}(\Lambda)=\sum_{\bm{v}\in \Lambda\bk \{\bm{0}\}}f(\bm{v}).$$

For any $\bm{\alpha}\in \R^n$ we denote by $Y_{\bm{\alpha}}$ the space of affine lattices of the form $(\Z^n+\bm{\alpha})g$ with $g\in G_n$. In this section we will prove a first moment and a second moment formula for the Siegel transform restricted to $Y_{\bm{\alpha}}$ with $\bm{\alpha}$ a rational vector. Our first observation is that similar to the space of unimodular lattices, when $\bm{\alpha}$ is rational the space $Y_{\bm{\alpha}}$ can also be parameterized by a certain homogeneous space $\G\bk G_n$ with $\G$ a conjugate of the congruence subgroup $\G_1^n(q)$. This way when $\bm{\alpha}$ is rational we can naturally endow $Y_{\bm{\alpha}}$ with the right $G_n$-invariant measure $\mu_n$ we fixed before.
\begin{Lem}\label{lem:id1}
Let $(\bm{p}, q)\in \Z^n\times \N$ with $\gcd(\bm{p},q)=1$, and let $\gamma_{\bm{p}}\in \G_n$ such that $\bm{p}=r\bm{e}_n\gamma_{\bm{p}}$ with $r=\gcd\left(\bm{p}\right)$, then $Y_{\frac{\bm{p}}{q}}$ can be parameterized by the homogeneous space $\gamma_{\bm{p}}^{-1}\G_1^n(q)\gamma_{\bm{p}}\bk G_n$  via the identification identifying $(\Z^n+\frac{\bm{p}}{q})g$ with $\gamma_{\bm{p}}^{-1}\G_1^n(q)\gamma_{\bm{p}}g$.
\end{Lem}
\begin{proof}
Consider the right multiplication action of $G_n$ on $Y_{\frac{\bm{p}}{q}}$ that $g\cdot \Lambda=\Lambda g^{-1}$ for $g\in G_n$ and $\Lambda\in Y_{\frac{\bm{p}}{q}}$. It is clear from the definition of $Y_{\frac{\bm{p}}{q}}$ that this action is transitive and it thus suffices to show that the stabilizer of $\Z^n+\frac{\bm{p}}{q}$ under this action is $\gamma_{\bm{p}}^{-1}\G_1^n(q)\gamma_{\bm{p}}$. The stabilizer clearly contains $\gamma_{\bm{p}}^{-1}\G_1^n(q)\gamma_{\bm{p}}$.  For the other containment, suppose $g\in G_n$ fixes $\Z^n+\frac{\bm{p}}{q}$, then $\Z^ng+\frac{\bm{p}}{q}g-\frac{\bm{p}}{q}=\Z^n$. This implies that $\Z^ng=\Z^n$ and $\frac{\bm{p}}{q}g-\frac{\bm{p}}{q}\in \Z^n$. The first condition implies that $g\in \G_n$ and the second condition implies that $\bm{p}g\equiv \bm{p} \Mod{q}$. Hence $r\bm{e}_n\gamma_{\bm{p}}g\equiv r\bm{e}_n\gamma_{\bm{p}} \Mod{q}$. By assumption $r=\gcd \left(\bm{p}\right)$ is coprime to $q$, thus we have  $\bm{e}_n\gamma_{\bm{p}}g\equiv \bm{e}_n\gamma_{\bm{p}} \Mod{q}$, or equivalently, $\gamma_{\bm{p}}g\gamma_{\bm{p}}^{-1}\in \G_1^n(q)$. Hence we have $g\in \gamma_{\bm{p}}^{-1}\G_1^n(q)\gamma_{\bm{p}}$ completing the proof.
\end{proof}
We now state our moment formulas for the Siegel transform restricted to the space $Y_{\frac{\bm{p}}{q}}$.

\begin{Thm}\label{thm:moment}
Let $n\geq 2$ be an integer. Let $(\bm{p}, q)\in \Z^n\times \N$ with $\gcd(\bm{p},q)=1$, and let $Y_{\frac{\bm{p}}{q}}$ be as above. Let $f$ be a bounded and compactly supported function on $\R^n$. Then we have
$$\int_{Y_{\frac{\bm{p}}{q}}}\hat{f}(\Lambda)d\mu_n(\Lambda)=\nu_{n,q}\vol(f).$$
If we further assume that $n\geq 3$, then
\begin{equation}\label{equ:secondsiegel}
\int_{Y_{\frac{\bm{p}}{q}}}|\hat{f}(\Lambda)|^2d\mu_n(\Lambda)=\nu_{n,q}|\vol(f)|^2+\frac{q^n}{\zeta(n)}\sum_{\substack{k_1\geq 1\\ \gcd(k_1,q)=1}}\sum_{\substack{k_2\in \Z\setminus\{0\} \\ k_2\equiv k_1\Mod{q}}}\vol(f_{k_1}\bar{f}_{k_2}),
\end{equation}
where $\nu_{n,q}=\frac{q^n\zeta_q(n)}{\zeta(n)}$ and for any $\lambda>0$, $f_{\lambda}(\bm{x}):=f(\lambda\bm{x})$ are as before.
\end{Thm}
\begin{rmk}
When $q=1,$ \thmref{thm:moment} gives an alternative proof to Rogers' second moment formula \cite{Rogers1955} on the space of unimodular lattices. We also note that since we will prove \thmref{thm:moment} using the moments formulas for incomplete Eisenstein series, we can apply \thmref{thm:moment} to indicator functions of any finite-volume subsets of $\R^n$, see the paragraph after \eqref{equ:choiceoftau}.
\end{rmk}

\subsection{Relating the Siegel transform with incomplete Eisenstein series}
Our strategy of proving \thmref{thm:moment} is to relate the Siegel transform restricted to $Y_{\frac{\bm{p}}{q}}$ to the incomplete Eisenstein series $\Theta_f^{n,q}$ and then apply the moment formulas \eqref{equ:firstincomp} and \eqref{equ:secondmoment} to these incomplete Eisenstein series. The following lemma rewrites the Siegel transform restricted to $Y_{\frac{\bm{p}}{q}}$ as an infinite sum of incomplete Eisenstein series.
\begin{Lem}\label{lem:relationsiegel}
Let $(\bm{p},q)\in \Z^n\times \N$ and $\gamma_{\bm{p}}\in \G_n$ be as in \lemref{lem:id1}. Let $f: \R^n\to \C$ be a bounded compactly supported function on $\R^n$. Then for any $\Lambda=(\Z^n+\frac{\bm{p}}{q})g\in Y_{\frac{\bm{p}}{q}}$ with $g\in G_n$ we have
$$\hat{f}(\Lambda)=\sum_{\substack{k\geq 1\\ \gcd(k, q)=1}}\Theta^{n,q}_{f_{\frac{k}{q}}}(\tau_{[k\bar{r}]}\gamma_{\bm{p}}g),$$
where for any $[\ell]\in (\Z/q\Z)^{\times}$, $\tau_{[\ell]}\in \G_n$ is defined as in \eqref{equ:tau}, and $f_{\frac{k}{q}}(\bm{x})=f(\frac{k\bm{x}}{q})$ as before.
\end{Lem}
\begin{proof}
First we note that for any $\Lambda=(\Z^n+\frac{\bm{p}}{q})g\in Y_{\frac{\bm{p}}{q}}$ we have
\begin{equation}\label{equ:prerelation}
\hat{f}(\Lambda)=\sum_{\bm{v}\in (\Z^n+\frac{\bm{p}}{q})\setminus \{\bm{0}\}}f(\bm{v}g)=\sum_{\bm{v}\in \cO(\bm{p})}f(\frac{\bm{v}g}{q}),
\end{equation}
where $\cO(\bm{p})=\left\{\bm{v}\in \Z^n\setminus \{\bm{0}\}\ |\ \bm{v}\equiv \bm{p} \Mod{q}\right\}$ is as before. Note that for any $\bm{v}\in \cO(\bm{p})$ there exists a unique $k\geq 1$ such that $\bm{v}=k\bm{v}'$ with $\bm{v}'\in \Z_{\rm{pr}}^n$. Thus $\gcd\left(k,q\right)=\gcd\left(k\bm{v}',q\right)=\gcd\left(\bm{p},q\right)=1$ implying that $\bm{v}'\in \cO(\bm{p};[k])$ with $\cO(\bm{p};[k])=\left\{\bm{w}\in \Z^n_{\rm pr}\ |\ \bm{w}\equiv \bar{[k]}\bm{p}\Mod{q}\right\}$ defined as before. Conversely, for any $k\geq 1$ coprime to $q$ and any $\bm{v}'\in \cO(\bm{p};[k])$ we have $k\bm{v}'\in \cO(\bm{p})$. We thus have the decomposition
$$\cO(\bm{p})=\bigsqcup_{[\ell]\in (\Z/q\Z)^{\times}}\bigsqcup_{\substack{k\geq 1 \\ [k]=[\ell]}}k\cO(\bm{p};[\ell]).$$
Moreover, since $\bm{p}=r\bm{e}_n\gamma_{\bm{p}}$, for any $\bm{w}\in \cO\left(\bm{p};[\ell]\right)$ we have $\bm{w}\gamma_{\bm{p}}^{-1}\equiv [r\bar{\ell}]\bm{e}_n\Mod{q}$ implying that
$$\cO(\bm{p};[\ell])=\cO(\bm{e}_n; [\ell \bar{r}])\gamma_{\bm{p}}=\cO(\bm{e}_n;[1])\tau_{[\ell \bar{r}]}\gamma_{\bm{p}},$$
where for the second identity we used the relation $\cO(\bm{e}_n;[\ell\bar{r}])=\cO(\bm{e}_n;[1])\tau_{[\ell\bar{r}]}$. Hence we get decomposition
$$\cO(\bm{p})=\bigsqcup_{[\ell]\in (\Z/q\Z)^{\times}}\bigsqcup_{\substack{k\geq 1 \\ [k]=[\ell]}}k\cO(\bm{e}_n;[1])\tau_{[\ell\bar{r}]}\gamma_{\bm{p}}.$$
Now use this decomposition and \eqref{equ:prerelation} to get
\begin{align*}
\hat{f}(\Lambda)
&=\sum_{[\ell]\in (\Z/q\Z)^{\times}}\sum_{\substack{k\geq 1\\ [k]= [\ell]}}\sum_{\bm{v}\in \cO(\bm{e}_n;[1])}f(\frac{k\bm{v}\tau_{[\ell\bar{r}]}\gamma_{\bm{p}}g}{q})\\
&=\sum_{[\ell]\in (\Z/q\Z)^{\times}}\sum_{\substack{k\geq 1\\ [k]= [\ell]}}\sum_{\gamma\in \G_1^n(q)\cap P_n\bk \G_1^n(q)}f(\frac{k\bm{e}_n\gamma\tau_{[\ell\bar{r}]}\gamma_{\bm{p}}g}{q})\\
&=\sum_{[\ell]\in (\Z/q\Z)^{\times}}\sum_{\substack{k\geq 1\\ [k]= [\ell]}}\Theta^{n,q}_{f_{\frac{k}{q}}}(\tau_{[\ell\bar{r}]}\gamma_{\bm{p}}g)=\sum_{\substack{k\geq 1\\ \gcd(k, q)=1}}\Theta^{n,q}_{f_{\frac{k}{q}}}(\tau_{[k\bar{r}]}\gamma_{\bm{p}}g),
\end{align*}
where for the second equality we used \lemref{lem:id2}.
\end{proof}

\begin{proof}[Proof of \thmref{thm:moment}]
Let $\cF_{\frac{\bm{p}}{q}}\subset G_n$ be a fundamental domain for $\gamma_{\bm{p}}^{-1}\G_1^n(q)\gamma_{\bm{p}}\bk G_n$ and note that $\cF_q:=\gamma_{\bm{p}}\cF_{\frac{\bm{p}}{q}}\subset G_n$ forms a fundamental domain for $X_{n,q}=\G_1^n(q)\bk G_n$. Then making a change of variable $\gamma_{\bm{p}}g\mapsto g$ we have
$$\int_{Y_{\frac{\bm{p}}{q}}}\hat{f}(\Lambda)d\mu_n(g)=\mathop{\sum_{k\geq 1}}_{\gcd(k,q)=1}\int_{\cF_q}\Theta^{n,q}_{f_{\frac{k}{q}}}(\tau_{[k\bar{r}]}g)d\mu_n(g).$$
Next, we note that for any $\gamma\in \G_n$, the shift $\gamma\cF_q$ forms a fundamental domain for $\gamma\G_1^n(q)\gamma^{-1}\bk G_n$. In particular, for any $[k]\in (\Z/q\Z)^{\times}$, since $\tau_{[k\bar{r}]}$ normalizes $\G_1^n(q)$, we get that $\tau_{[k\bar{r}]}\cF_q$ is also a fundamental domain for $X_{n,q}$. Thus for each $[k]\in (\Z/q\Z)^{\times}$ making a change of variable $\tau_{[k\bar{r}]}g\mapsto g$ and applying \eqref{equ:firstincomp} and using \lemref{lem:index} we get
$$\int_{Y_{\frac{\bm{p}}{q}}}\hat{f}(\Lambda)d\mu_n(g)=\sum_{\substack{k\geq 1\\ \gcd(k,q)=1}}\frac{\vol(f_{\frac{k}{q}})}{\zeta(n)}=\frac{q^n\zeta_q(n)}{\zeta(n)}\vol(f)=\nu_{n,q}\vol(f)$$
as claimed.
%

Next for the second moment, making a change of variable $\gamma_{\bm{p}}g\mapsto g$ we get
$$\int_{Y_{\frac{\bm{p}}{q}}}|\hat{f}(\Lambda)|^2d\mu_n(\Lambda)=\sum_{\substack{k_i\geq 1\\ \gcd(k_i,q)=1\\ i=1,2}}\int_{\cF_q}\Theta^{n,q}_{f_{\frac{k_1}{q}}}(\tau_{[k_1\bar{r}]}g)\Theta^{n,q}_{\bar{f}_{\frac{k_2}{q}}}(\tau_{[k_2\bar{r}]}g)d\mu_n(g).$$
Making a change of variable $\tau_{[k_1\bar{r}]}g\mapsto g$, noting that $\tau_{[k_1\bar{r}]}\cF_q$ is also a fundamental domain for $X_{n,q}$ and using the relation \eqref{equ:choiceoftau}, we have $\int_{Y_{\frac{\bm{p}}{q}}}|\hat{f}(\Lambda)|^2d\mu_n(\Lambda)$ equals
\begin{align*}
&\sum_{\substack{k_i\geq 1\\ \gcd(k_i,q)=1\\ i=1,2}}\int_{\cF_q}\Theta^{n,q}_{f_{\frac{k_1}{q}}}(g)\Theta^{n,q}_{\bar{f}_{\frac{k_2}{q}}}(\tau_{[k_1\bar{r}]}^{-1}\tau_{[k_2\overline{r}]}g)d\mu_n(g)
=\sum_{\substack{k_i\geq 1\\ \gcd(k_i,q)=1\\ i=1,2}}\int_{\cF_q}\Theta^{n,q}_{f_{\frac{k_1}{q}}}(g)\Theta^{n,q}_{\bar{f}_{\frac{k_2}{q}}}(\tau_{[\bar{k}_1k_2]}g)d\mu_n(g).
\end{align*}
Now applying \thmref{thm:secondmomentinein} to the functions $f_{\frac{k_1}{q}}$ and $\bar{f}_{\frac{k_2}{q}}$ we get that $\int_{Y_{\frac{\bm{p}}{q}}}|\hat{f}(\Lambda)|^2d\mu_n(\Lambda)$ equals
\begin{align*}
&\sum_{\substack{k_i\geq 1\\ \gcd(k_i,q)=1\\ i=1,2}}\left(\frac{\vol\left(f_{\frac{k_1}{q}}\right)\vol\left(\bar{f}_{\frac{k_2}{q}}\right)}{q^n\zeta(n)\zeta_q(n)} +\frac{1}{\zeta(n)}\left(\delta_{[k_1][k_2]}\vol\left(f_{\frac{k_1}{q}}\bar{f}_{\frac{k_2}{q}}\right)+\delta_{[k_1][-k_2]}\vol\left(f_{\frac{k_1}{q}}\tilde{\bar{f}}_{\frac{k_2}{q}}\right)\right)\right)\\
&= \frac{|\vol(f)|^2}{q^n\zeta(n)\zeta_q(n)}\sum_{\substack{k_i\geq 1\\ \gcd(k_i,q)=1\\ i=1,2}}\frac{q^{2n}}{k_1^nk_2^n}+\frac{q^n}{\zeta(n)}\sum_{\substack{k_1\geq 1\\ \gcd(k_1,q)=1}}\sum_{\substack{k_2\in \Z\setminus\{0\} \\ k_2\equiv k_1\Mod{q}}}\vol(f_{k_1}\bar{f}_{k_2})\\
&=\nu_{n,q}|\vol(f)|^2+\frac{q^n}{\zeta(n)}\sum_{\substack{k_1\geq 1\\ \gcd(k_1,q)=1}}\sum_{\substack{k_2\in \Z\setminus\{0\} \\ k_2\equiv k_1\Mod{q}}}\vol(f_{k_1}\bar{f}_{k_2}),
\end{align*}
finishing the proof, where for the first equality we used that $\tilde{\bar{f}}_{k_2/q}=\bar{f}_{-k_2/q}$ and did a change of variable $-k_2\mapsto k_2$ and for the second equality we used the identity that
\begin{displaymath}
\frac{1}{q^n\zeta(n)\zeta_q(n)}\sum_{\substack{k_i\geq 1\\ \gcd(k_i,q)=1\\ i=1,2}}\frac{q^{2n}}{k_1^nk_2^n}=\frac{q^n\zeta_q(n)}{\zeta(n)}=\nu_{n,q}.\qedhere
\end{displaymath}
\end{proof}
\subsection{An alternative proof to the second moment formula}
As illustrated to us by an anonymous referee, the second moment formula in \thmref{thm:moment} can also be proved differently using results from \cite[Section 7]{MarklofStrombergsson2010}. Following the referee's suggestions we sketch this alternative proof below.
We first need to introduce some notations and definitions from \cite{MarklofStrombergsson2010}.
We will fix $\bm{\alpha}=\frac{\bm{p}}{q}$. Let 
$$\G(q):=\left\{\gamma\in \G_n\ |\ \gamma\equiv I_n\Mod{q}\right\}$$ 
denote the principle congruence subgroup of level $q$ and let $X_q=\G(q)\bk G_n$. For any $\bm{y}\in \R^n\setminus\{\bm{0}\}$ let 
$$X_q(\bm y )=\{\G(q) g\in X_q\ |\ \bm{y}\in (\Z^n+\alpha)g\}. $$
It was shown in  \cite{MarklofStrombergsson2010}  that $X_q(\bm y)$ is an embedded submanifold of $X_q$ and that it carries an invariant Borel measure $\nu_{\bm y}$ (see the paragraph before \cite[Lemma 7.2]{MarklofStrombergsson2010} for the precise definition).  We will now show how the second moment formula follows from the following two properties of this measure.

The first is \cite[Proposition 7.6]{MarklofStrombergsson2010} stating that for any compactly supported bounded measurable function $f:\R^n\to \R$ and for any $\bm{y}\in \R^n\setminus\{\bm{0}\},$ 
\begin{equation}\label{e:p7.6}
\int_{X_q(\bm y)} \hat{f}((\Z^n+\bm \alpha)g)d\nu_{\bm y}(g)=\vol(f)+\mathop{\sum_{t\geq 1}}_{\gcd(t,q)=1}t^{-n}\mathop{\sum_{a\in   q\Z+t}}_{\gcd(a,t)=1} f\left(\frac{a}{t}\bm y\right),
\end{equation}
and the second is  \cite[(7.25)]{MarklofStrombergsson2010} stating that for any Borel set $\cE\subset X_q$ and for any Borel subset $U\subset \R^n\setminus \{\bm{0}\},$
$$\int_U\nu_{\bm{y}}\left(\cE\cap X_q(\bm{y})\right)d\bm{y}=\sum_{\bm{v}\in (\Z^n+\bm{\alpha})\setminus\{\bm{0}\}}\widetilde{\mu}_n(\cF\cap \cE_{\bm{v}}),$$
where $\widetilde{\mu}_n=\mu_n/[\G_n, \G(q)]$ is the normalized Haar measure on $X_q$, $\cF\subset G_n$ is a fundamental domain for $X_q$ and for any $\bm{v}\in (\Z^n+\bm{\alpha})\setminus \{\bm{0}\}$
$$\cE_{\bm{v}}:=\{g\in G_n\ |\ \G(q)g\in \cE,\ \bm{v}g\in U\}.$$

Using the definition of $\cE_{\bm{v}}$ and changing the order of summation and integration we have the following identity
\begin{equation}\label{equ:prefosi}
\int_{X_q}\sum_{\bm{v}\in (\Z^n+\bm{\alpha})\setminus \{\bm{0}\}}\chi_{U}(\bm{v}g)\chi_{\cE}(\G(q)g)d\widetilde{\mu}_n(g)=\int_U\nu_{\bm{y}}\left(\cE\cap X_q(\bm{y})\right)d\bm{y},
\end{equation}
where both sides of the equality are allowed to be $\infty$. Using a standard approximation argument from measure theory, one can deduce from \eqref{equ:prefosi} the following equality for arbitrary Borel measurable functions $F: X_q\times (\R^n\setminus\{\bm{0}\})\to \R_{\geq 0}$, that is
\begin{equation}\label{equ:presieg}
\int_{X_q}\sum_{\bm{v}\in (\Z^n+\bm{\alpha})\setminus \{\bm{0}\}}F(\G(q)g, \bm{v}g)d\widetilde{\mu}_n(g)=\int_{\R^n\setminus \{\bm{0}\}}\int_{X_q(\bm{y})}F(\G(q)g, \bm{y})d\nu_{\bm{y}}(g)d\bm{y}.
\end{equation}

Now, we first assume that $f$ is non-negative. Using the identification between $Y_{\bm{\alpha}}$ and the homogeneous space $\gamma_{\bm{p}}^{-1}\G_1^n(q)\gamma_{\bm{p}}\bk G_n$ and the facts that $\G(q)<\G_1^n(q)$ and that $\G(q)$ is a normal subgroup of $\G_n$, we have that $\hat{f}\big|_{Y_{\bm{\alpha}}}$ is left $\G(q)$-invariant, so we can think of $\hat{f}\big|_{Y_{\bm{\alpha}}}$ as a function on $X_q$. Now unfolding one of the factors and  applying \eqref{equ:presieg} with $F(\G(q)g, \bm{y})=\hat{f}(\G(q)g)f(\bm{y})$ we get that 
\begin{align*}
\frac{1}{\nu_{n,q}}\int_{Y_{\bm{\alpha}}}\left|\hat{f}(\Lambda)\right|^2d\mu_n(\Lambda)&=\int_{X_q}\sum_{\bm{v}\in (\Z^n+\bm{\alpha})\setminus \{\bm{0}\}}f(\bm{v}g)\hat{f}(\G(q)g)d\widetilde{\mu}_n(g)\\
&=\int_{\R^n\setminus \{\bm{0}\}}\int_{X_q(\bm{y})}\hat{f}(\G(q)g)f(\bm{y})d\nu_{\bm{y}}(g)d\bm{y}\\
&=\int_{\R^n\setminus \{\bm{0}\}}f(\bm y) \left(\int_{X_q(\bm{y})}\hat{f}(\G(q)g)d\nu_{\bm{y}}(g)\right)d\bm{y}.
\end{align*}
Next, note that $\hat{f}(\G(q)g)=\hat{f}((\Z^n+\bm \alpha)g)$ and apply \eqref{e:p7.6} to the inner integral to get that

\begin{align*}
\frac{1}{\nu_{n,q}}\int_{Y_{\bm{\alpha}}}\left|\hat{f}(\Lambda)\right|^2d\mu_n(\Lambda)&=
\int_{\R^n\setminus\{\bm{0}\}}f(\bm y) \left(\int_{\R^n}f(\bm{x})d\bm{x}+\sum_{\substack{t\geq 1\\ \gcd(t,q)=1}}t^{-n}\sum_{\substack{a\in (t+q\Z)\setminus\{0\}\\ \gcd(a,t)=1}}f(\frac{a}{t}\bm{y})\right)d\bm{y}\\
&=\vol(f)^2+\sum_{\substack{t\geq 1\\ \gcd(t,q)=1}}\sum_{\substack{a\in (t+q\Z)\setminus\{0\}\\ \gcd(a,t)=1}}\int_{\R^n}f(a\bm{y})f(t\bm{y})d\bm{y}\\
&=\vol(f)^2+\frac{1}{\zeta_q(n)}\sum_{\substack{t\geq 1\\ \gcd(t,q)=1}}\sum_{\substack{a\in (t+q\Z)\setminus\{0\}\\ \gcd(a,t)=1}}\sum_{\substack{\delta\geq 1\\ \gcd(\delta, q)=1}}\int_{\R^n}f(\delta a\bm{y})f(\delta t\bm{y})d\bm{y},
\end{align*}
where for the second equality we did a change of variable $\bm{y}/t\mapsto \bm{y}$ and for the last equality we used the identity $\zeta_q(n)=\sum_{\substack{\delta\geq 1\\ \gcd(\delta, q)=1}}\delta^{-n}$ and did a change of variable $\bm{y}\mapsto \delta \bm{y}$. 

Finally, making  a change of variables $\delta a\mapsto k_1$ and $\delta t\mapsto k_2$ gives the second moment formula \eqref{equ:secondsiegel}  concluding the proof for $f$ a non-negative function. For a general complex-valued function $f$, one can check that the right hand side (and hence also the left hand side) of \eqref{equ:secondsiegel}  absolutely converges when $f$ is bounded and compactly supported. Thus by splitting $f$ into real and imaginary parts and further into positive and negative parts, one can extend \eqref{equ:secondsiegel} for a complex-valued function $f$.

\subsection{Application to discrepancies}
As a direct consequence of our moment formulas, we have the following mean square bound for the discrepancy function. For any finite-volume Borel set $A\subset \R^n$ and for any affine lattice $\Lambda\in Y_{\frac{\bm{p}}{q}}$ we define the discrepancy function as 
$$D(\Lambda,A):=\left|\#(\Lambda\cap A)-\vol(A)\right|.$$
\begin{Cor}\label{cor:discrepancybound}
Keep the notation as above. For any finite-volume Borel set $A\subset \R^n$ with $\vol(A)>1$ we have
$$\int_{Y_{\frac{\bm{p}}{q}}}\left|D(\Lambda,A)\right|^2d\mu_n(\Lambda)\ll_n q^n\vol(A),$$
where the bounding constant only depends on $n$. 
\end{Cor}
\begin{Rmk}
We note that the assumption that $\vol(A)>1$ is only needed to handle the case when $q=1$ and $\bm{0}\in A$.
\end{Rmk}
\begin{proof}[Proof of \corref{cor:discrepancybound}]
Using the first moment formula in \thmref{thm:moment}, the assumption that $\vol(A)>1$ and the relation that for any $\Lambda\in Y_{\frac{\bm{p}}{q}}$ 
$$0\leq \#\left(\Lambda\cap A\right)-\#\left(\Lambda\cap \left(A\setminus\{\bm{0}\}\right)\right)\leq 1,$$ 
we may assume without loss of generality that $\bm{0}\notin A$. Let $f=\chi_A$ be the indicator function of $A$, and for any $k\geq 1$ let $f_k(\bm{x})=f(k\bm{x})$ as before. We note that since $\bm{0}\notin A$, for any $\Lambda\in Y_{\frac{\bm{p}}{q}}$ we have $\hat{f}(\Lambda)=\hat{\chi}_A(\Lambda)=\#(\Lambda\cap A)$. Hence integrating $D(\Lambda, A)$ over $Y_{\frac{\bm{p}}{q}}$ we get
\begin{align*}
\int_{Y_{\frac{\bm{p}}{q}}}D(\Lambda,A)^2d\mu_n(\Lambda)&
=\int_{Y_{\frac{\bm{p}}{q}}}\left|\hat{f}(\Lambda)\right|^2d\mu_n(\Lambda)-2\vol(A)\int_{Y_{\frac{\bm{p}}{q}}}\hat{f}(\Lambda)d\mu_n(\Lambda)+\vol(A)^2\mu_n(Y_{\frac{\bm{p}}{q}})\\
&=\int_{Y_{\frac{\bm{p}}{q}}}\left|\hat{f}(\Lambda)\right|^2d\mu_n(\Lambda)-\nu_{n,q}\vol(A)^2\\
&=\frac{q^n}{\zeta(n)}\sum_{\substack{k_1\geq 1\\ \gcd(k_1,q)=1}}\sum_{\substack{k_2\in \Z\setminus\{0\} \\ k_2\equiv k_1\Mod{q}}}\vol(f_{k_1}f_{k_2})\\
&\leq \frac{q^n}{\zeta(n)}\sum_{\substack{k_1\geq 1\\ \gcd(k_1,q)=1}}\sum_{\substack{k_2\in \Z\setminus\{0\} \\ k_2\equiv k_1\Mod{q}}}\frac{\vol(A)}{k_1^{\frac{n}{2}}k_2^{\frac{n}{2}}}\leq\frac{2\zeta(\frac{n}{2})^2}{\zeta(n)}q^n\vol(A),
\end{align*}
where for the second equality we applied the first moment formula  
and \lemref{lem:id1}, for the third equality we applied the second moment formula in \thmref{thm:moment}, for the first inequality we applied the Cauchy-Schwarz inequality, and for the last inequality we used the following bound 
\begin{displaymath}
\sum_{\substack{k_1\geq 1\\ \gcd(k_1,q)=1}}\sum_{\substack{k_2\in \Z\setminus\{0\} \\ k_2\equiv k_1\Mod{q}}}\frac{1}{k_1^{\frac{n}{2}}k_2^{\frac{n}{2}}}\leq 2\sum_{k_1,k_2\geq 1}\frac{1}{k_1^{\frac{n}{2}}k_2^{\frac{n}{2}}}=2\zeta(\frac{n}{2})^2<\infty.\qedhere
\end{displaymath}
\end{proof}



\section{Rational shifts and generic forms}\label{sec:rational}
In this section we give the proof of  \thmref{thm:rationalshift}. First we note that in view of \rmkref{rmk:scale} it suffices to consider quadratic forms with unit determinant. Next we note that with our new moment formulas, the proof is almost identical to the ones given in \cite{KelmerYu2018b}. Here we only collect the necessary ingredients needed for the proof and refer the reader to \cite[Theorem 6]{KelmerYu2018b} for the details. 
Let us first fix some notation. For any $n=p_1+p_2\geq 3$ with $p_1,p_2\geq 1$ we denote by $\cQ_{p_1,p_2}$ the space of unit determinant quadratic forms of signature $(p_1,p_2)$. We fix a base quadratic form $Q_0\in\cQ_{p_1,p_2}$ such that
$$Q_0(\bm{v})=\sum_{i=1}^{p_1}v_i^2-\sum_{i=p_1+1}^nv_i^2.$$
As mentioned in the introduction, any other quadratic form $Q$ in $\cQ_{p_1,p_2}$ can be written in the form $Q(\bm{v})=g\cdot Q_0(\bm{v}):=Q_0(\bm{v}g)$ for some $g\in G_n$.

As in \cite{KelmerYu2018b}, we will prove  \thmref{thm:rationalshift} by studying a lattice point counting problem. More precisely, for any inhomogeneous form $Q_{\bm{\alpha}}$ with $Q=g\cdot Q_0\in \cQ_{p_1,p_2}$ and $\bm{\alpha}\in\R^n$, for any target set $I\subset \R$ and for any $t>0$, we have
\begin{equation}\label{equ:countingrelation}
\cN_{Q_{\bm{\alpha}},I}(t)=\#\left(\Z^n\medcap Q_{\bm{\alpha}}^{-1}(I)\medcap B_t\right)=\left((\Z^n+\bm{\alpha})g\medcap Q_0^{-1}(I)\medcap (B_t+\bm{\alpha})g\right),
\end{equation}
where $\cN_{Q_{\bm{\alpha}},I}(t)$ is the counting function defined as in \eqref{equ:countingfunction}, $B_t\subset \R^n$ is the open Euclidean ball as fixed in the introduction, and for the second equality we used that 
$$Q_{\bm{\alpha}}^{-1}(I)=Q_0^{-1}(I)g^{-1}-\bm{\alpha}.$$ 
To simplify notation for fixed $\bm{\alpha}\in\R^n$ let us define for any $g\in G_n$, for any subset $I\subset \R$ and for any $t>0$
$$A_{g,I,t}:=Q_0^{-1}(I)\medcap (B_t+\bm{\alpha})g.$$ 
In view of the above relation, we are thus interested in proving an asymptotic power saving bound for the discrepancy function $D\left((\Z^n+\bm{\alpha})g, A_{g,I_t,t}\right)$ with $\bm{\alpha}$ and $I_t$ as in  \thmref{thm:rationalshift}. We note that as $g$ runs through $G_n$, the set of affine lattices we consider here is exactly the space $Y_{\bm{\alpha}}$ as defined in \secref{sec:Siegel}. In particular, when $\bm{\alpha}$ is rational 
we have the following bound for measures of the sets of $\Lambda\in Y_{\bm{\alpha}}$ with large discrepancies. We refer the reader to \cite[Lemma 2.2]{KelmerYu2018b} for the proof of the case when $q=1$, and we note that the same proof carries over for general $q\in \N$ after replacing \cite[(2.4)]{KelmerYu2018b} with our mean square bound (\corref{cor:discrepancybound}) and noting that the number of copies of fundamental domains of $X_{n,q}$ needed to cover a compact set $\mathcal{K}\subset G_n$ is bounded from above by that of $X_{n,1}$.
\begin{Lem}\label{lemma:MAT}
Let $\bm{\alpha}=\frac{\bm{p}}{q}\in \Q^n$ be a fixed rational vector with $(\bm{p},q)\in \Z^n\times \N$ and $\gcd\left(\bm{p},q\right)=1$. Fix $\cK\subset G_n$ a compact subset with positive measure. For any bounded and measurable subset $A\subset \R^n$ with $\vol(A)>1$ and any $T>0$ 
\begin{equation}\label{equ:discrepancy}
\mu_n\left(\cM^{(\cK,\bm{\alpha})}_{A,T}\right)\ll_{\cK,n} \frac{q^n\vol(A)}{T^2},
\end{equation}
where $\cM^{(\cK,\bm{\alpha})}_{A,T}:=\left\{g\in \cK\ \left|\ D\left((\Z^n+\bm{\alpha})g,A\right)\geq T\right.\right\}$.
\end{Lem}

Another key ingredient in our proof is an effective volume estimate for sets of the form $Q^{-1}(I)\cap B_t$ which holds for all quadratic forms $Q\in \cQ_{p_1,p_2}$. For this we record the following effective volume estimate from \cite[Theorem 5]{KelmerYu2018b}. We note that the volume estimate in \cite{KelmerYu2018b} holds for a more general family of homogeneous polynomials of a fixed even degree. In particular, taking the degree to be two we get the following.
\begin{Thm}[{{\cite[Theorem 5]{KelmerYu2018b}}}]\label{thm:volumeestimate}
Let $n=p_1+p_2\geq 3$ with $p_1, p_2\geq 1$, and let $N\geq 1$. For $Q\in \cQ_{p_1,p_2}$ and $I\subseteq [-N,N]$ measurable, there exists $c_Q>0$ such that  for any $t> 2N^{1/2}$ we have 
$$\vol(Q^{-1}(I)\cap B_t)=c_Q |I|t^{n-2}+O_Q(|I| N^{1/2}t^{n-3}\log(t)),$$
where the implied constant is uniform over compact sets and the $\log(t)$ factor is only needed when $n=3$.
\end{Thm}
For our application the interval $I$ is shrinking to a fixed point so we may ignore the dependence on $N$. Moreover, using the estimates 
$B_{t-\|\bm{\alpha}\|}\subset B_t+\bm{\alpha}\subset B_{t+\|\bm{\alpha}\|}$ and $(t\pm \|\bm{\alpha}\|)^a=t^a+O_{\bm{\alpha}}(t^{a-1})$ for any $t>2\|\bm{\alpha}\|$ and $a\in\R$, 
we get the following.
\begin{Cor}\label{cor:volumeestimate}
Keep the notation as in \thmref{thm:volumeestimate}. For $A_{g,I,t}:=Q_0^{-1}(I)\cap (B_t+\bm{\alpha})g$ as above with $I\subset [-N,N]$ we have for $t>2N^{1/2}$
$$\vol(A_{g,I,t})=c_Q |I|t^{n-2}+O_{g,\bm{\alpha}}(|I| t^{n-3}\log(t)),$$
where $Q=g\cdot Q_0$ and the $\log(t)$ term is needed only when $n=3$.
\end{Cor}
Following the same arguments as in the proof of \cite[Theorem 6]{KelmerYu2018b} with slight modifications replacing \cite[Lemma 2.2]{KelmerYu2018b} by \lemref{lemma:MAT}, 
and using the aforementioned two estimates,
gives the following bound for the discrepancies $D\left(\left(\Z^n+\bm{\alpha}\right)g,A_{g,I_t,t}\right)$. 
\begin{Thm}\label{thm:latticecounting}
Keep the assumptions as in \thmref{thm:rationalshift} and keep the notation as above. Then there exists some $\delta\in (0,1)$ such that for $\mu_n$-a.e. $g\in G_n$ there exists $t_{g,\bm{\alpha}}>0$ such that for all $t\geq t_{g,\bm{\alpha}}$
$$D\left(\left(\Z^n+\bm{\alpha}\right)g,A_{g,I_t,t}\right)<\vol(A_{g,I_t,t})^{\delta}.$$
\end{Thm}

The proof of  \thmref{thm:rationalshift} and \corref{c:cong} now easily follows.
\begin{proof}[Proof of \thmref{thm:rationalshift}]
Let $\bm{\alpha}\in \Q^n$ be fixed. Using the scaling trick in \rmkref{rmk:scale}, it suffices to prove this theorem for generic unit determinant forms. Fix $p_1, p_2\geq 1$ with $p_1+p_2=n$ and let $Q_0$ be the fixed unit determinant form of signature $(p_1,p_2)$ as above. Then by \thmref{thm:latticecounting} and the relation  \eqref{equ:countingrelation} we get that there is some $\delta\in(0,1)$ such that for $\mu_n$-almost all $g\in G_n$ there is $t_{g,\bm\alpha}>0$ such that for $Q=g\cdot Q_0$ and for all $t\geq t_{g,\bm \alpha}$
$$\cN_{Q_{\bm \alpha}, I_t}(t)=\#\left((\Z^n+\bm{\alpha})g\cap A_{g,I_t,t}\right)=\vol(A_{g,I_t,t})+O(\vol(A_{g,I_t,t})^\delta).$$
Thus by \corref{cor:volumeestimate}
$$\vol(A_{g,I_t,t})=c_Q |I_t|t^{n-2}+O_Q(|I_t| t^{n-3}\log(t)))$$
we get that for almost all unit determinant non-degenerate quadratic forms $Q$ of signature $(p_1,p_2)$,
$$\cN_{Q_{\bm\alpha}, I_t}(t)=c_Q |I_t|t^{n-2}+O_{Q,\bm{\alpha}}(t^{n-2-\kappa-\nu}),$$
with some positive $\nu<\min\{(n-2-\kappa)(1-\delta),1\}$. Since this holds for any signature $(p_1,p_2)$ with $p_1, p_2\geq 1$, the result holds for almost all unit determinant non-degenerate indefinite forms in $n$ variables.
\end{proof}

\begin{proof}[Proof of \corref{c:cong}]
For any $q\in \N$ and $\bm p\in (\Z/q\Z)^n$ let $\bm\alpha=\frac{\bm p}{q}$. For any $\bm{v}\in \Z^n$ write $\bm{w}=q\bm{v}+\bm{p}$ and use the observation \eqref{equ:obs} to get that for any $I\subset\R$, any $t>0$ and any quadratic form $Q$
$$\#\{\bm w\in \Z^n\ |\ Q(\bm w)\in I, \bm w\equiv\bm p\Mod{q}, \|\bm w\|\leq t\}=\#\{\bm v\in \Z^n\ |\ Q_{\bm\alpha}(\bm v)\in q^{-2}I, \|\bm v+\bm \alpha\|\leq q^{-1}t\}.$$
Now, the left hand side is bounded from above and below respectively by $\cN_{Q_{\bm \alpha},q^{-2}I}\left(\frac{t\pm \sqrt{n}}{q}\right)$, so the result follows by applying the power saving formula \eqref{equ:countingresult} to $\cN_{Q_{\bm \alpha},q^{-2}I_t}\left(\frac{t\pm \sqrt{n}}{q}\right)$ for the full measure set of $Q$ satisfying \eqref{equ:countingresult}.
\end{proof}


\section{Irrational shifts and generic forms}
In this section we will prove \thmref{thm:irrational} by proving a slightly more general result (\thmref{thm:moregeneral}). The proof uses some properties of simultaneous Diophantine approximations, which we now review. 
\subsection{Backgrounds on simultaneous Diophantine approximation}\label{sec:dio}
For any $\bm{\alpha}\in \R^n$ we denote by 
$$\langle \bm{\alpha}\rangle:=\inf_{\bm{v}\in \Z^n}\|\bm{\alpha}-\bm{v}\|_{\infty}$$ 
to be the shortest distance of $\bm{\alpha}$ to integer points with respect to the supremum norm $\|\cdot\|_{\infty}$ on $\R^n$. Recall that for any irrational $\bm{\alpha}\in\R^n$ the \textit{Diophantine exponent $\omega_{\bm{\alpha}}$} (respectively \textit{uniform Diophantine exponent $\hat{\omega}_{\bm{\alpha}}$}) of $\bm{\alpha}$ is the supremum of $\nu>0$ for which the system of inequalities
\begin{equation}\label{equ:dioexp}
\langle q\bm{\alpha}\rangle < t^{-\nu}\quad \textrm{and}\quad |q|< t
\end{equation}
has nontrivial integer solutions in $q$ for an unbounded set of $t>0$ (respectively for all sufficiently large $t>0$). It is clear that $\omega_{\bm{\alpha}}\geq \hat{\omega}_{\bm{\alpha}}$ and by the generalized Dirichlet's theorem on simultaneous Diophantine approximation (see e.g. \cite[\S \rom{1}.5]{Cassels1957}) we have that $\hat{\omega}_{\bm{\alpha}}\geq \frac{1}{n}$ for any irrational $\bm{\alpha}\in\R^n$. Moreover, by the Borel-Cantelli lemma we have that $\omega_{\bm{\alpha}}= \hat{\omega}_{\bm{\alpha}}= \frac{1}{n}$ for Lebesgue almost every $\bm{\alpha}\in\R^n$. We call the vectors $\bm{\alpha}$ for which  $\omega_{\bm{\alpha}}<\infty$, \emph{Diophantine vectors}. 

For any irrational $\bm{\alpha}=(\alpha_1,\ldots,\alpha_n)\in\R^n$ let $\textrm{Span}_{\Q}(\bm{\alpha})$ be the linear span of $1,\alpha_1,\ldots, \alpha_n$ over the rationals, and we denote by $d_{\bm{\alpha}}=\dim \textrm{Span}_{\Q}(\bm{\alpha})-1$. It is clear that $d_{\bm{\alpha}}\leq n$ and we note that $d_{\bm{\alpha}}=n$ if and only if $1,\alpha_1,\ldots, \alpha_n$ are linearly independent over $\Q$. We call $\bm{\alpha}\in\R^n$  \textit{totally irrational} if $d_{\bm{\alpha}}=n$. Choose a basis $\{1,\lambda_1,\ldots,\lambda_{d_{\bm{\alpha}}}\}$ for $\textrm{Span}_{\Q}(\bm{\alpha})$ and let $\bm{\lambda}=(\lambda_1,\ldots,\lambda_{d_{\bm{\alpha}}})\in \R^{d_{\bm{\alpha}}}$. It is then not hard to see that $\omega_{\bm{\alpha}}=\omega_{\bm{\lambda}}$ and $\hat{\omega}_{\bm{\alpha}}=\hat{\omega}_{\bm{\lambda}}$, implying that $\omega_{\bm{\alpha}}\geq \hat{\omega}_{\bm{\alpha}}\geq \frac{1}{d_{\bm{\alpha}}}$.
\subsubsection{Best simultaneous Diophantine approximation denominators}
Following \cite{Chevallier2013}, for any $\bm{\alpha}\in \R^n$ we say a positive integer $q$ is a \textit{best simultaneous Diophantine approximation denominator} of $\bm{\alpha}$ if $\langle q\bm{\alpha}\rangle < \langle l\bm{\alpha}\rangle$ for any $1\leq l< q$. For any irrational $\bm{\alpha}\in \R^n$, the set of best simultaneous Diophantine approximation denominators is infinite, and we thus get  an increasing sequence $q_0=q_0(\bm{\alpha})=1<q_1=q_1(\bm{\alpha})<\ldots <q_k=q_k(\bm{\alpha})<\ldots$ of best simultaneous Diophantine approximation denominators of $\bm{\alpha}$. Throughout this section, we will assume $\bm{\alpha}\in \R^n$ to be irrational. For each $q_k$, there exists a unique $\bm{p}_k\in \Z^n$ such that $\langle q_k\bm{\alpha}\rangle =\|q_k\bm{\alpha}-\bm{p}_k\|_{\infty}$ and note that it is clear from minimality that $\gcd (\bm{p}_k,q_k)=1$. For later use, we note that Lagarias \cite{Lagarias1982} showed that $\{q_k\}_{k\in\N}$ satisfies the relations $q_{k+2^n}\geq q_{k+1}+q_k$ and thus there exists some $c>1$ (depending on $\bm{\alpha}$) such that $q_k\geq c^k$. For instance, we can take $c>1$ sufficiently small such that $c\leq \min\left\{ q_i^{1/i}\ |\ 1\leq i\leq 2^n\right\}$ and $c+1\geq c^{2^n}$ and prove the result by induction. We call $\bm{r}_k:=\frac{\bm{p}_k}{q_k}$ the $k$th \textit{partial convergents} of $\bm{\alpha}$. Note that the aforementioned generalized Dirichlet's theorem implies that
\begin{equation}\label{equ:dirichlet}
\|\bm{\alpha}-\bm{r}_k\|_{\infty}<\frac{1}{q_kq_{k+1}^{\frac{1}{n}}}.
\end{equation}
The following simple lemma gives a relation between the Diophantine exponent $\bm{\alpha}$ and its sequence of best simultaneous Diophantine approximation denominators $\{q_k\}_{k\in\N}$.
\begin{Lem}\label{lem:badlyapproximable}
Keep the notation as above. Let $\bm{\alpha}\in \R^n$ be a Diophantine vector. Then for any $\nu>\omega_{\bm{\alpha}}$ the sequence $\{\frac{q_{k+1}}{q_k^{n\nu}}\}_{k\in \N}$ is bounded with the bounding constant depending on $\bm{\alpha}$ and $\nu$.
\end{Lem}
\begin{proof}
First we note that $\omega_{\bm{\alpha}}$ is the supremum of $\nu>0$ for which the inequality $\langle q\bm{\alpha}\rangle <q^{-\nu}$
has infinitely many integer solutions in $q$. Thus
$$\omega_{\bm{\alpha}}=\sup\left\{\nu>0\ |\ \liminf\limits_{q\to\infty}q^{\nu}\langle q\bm{\alpha}\rangle<\infty\right\}=\inf\left\{\nu>0\ |\ \liminf\limits_{q\to\infty}q^{\nu}\langle q\bm{\alpha}\rangle=\infty\right\},$$
where for the second equality we used the fact that $\liminf\limits_{q\to\infty}q^{\nu}\langle q\bm{\alpha}\rangle=0$ for any $\nu<\omega_{\bm{\alpha}}$ and $\liminf\limits_{q\to\infty}q^{\nu}\langle q\bm{\alpha}\rangle=\infty$ for any $\nu>\omega_{\bm{\alpha}}$.
Thus for any $\nu>\omega_{\bm{\alpha}}$ there exists $\varepsilon>0$ such that $\inf_{q\in \N}q^{\nu}\langle q\bm{\alpha} \rangle>\varepsilon$. In particular, taking $q=q_k$ and using \eqref{equ:dirichlet} we have $\varepsilon<q_k^{\nu}\langle q_k\bm{\alpha} \rangle<\frac{q^{\nu}_k}{q^{\frac{1}{n}}_{k+1}}$ implying that $\frac{q_{k+1}}{q^{n\nu}_k}<\varepsilon^{-n}$ is bounded. 
\end{proof}
The following lemma gives an interpretation of the uniform Diophantine exponent $\hat{\omega}_{\bm{\alpha}}$ in terms of the sequence $\{q_k\}_{k\in\N}$. We note that such an interpretation was used in \cite[Lemma 2.1]{KleinbockWadleigh2018}.
\begin{Lem}\label{lem:bestapp}
Keep the notation as above and let $\bm{\alpha}\in\R^n$ be irrational. Then we have
$$\hat{\omega}_{\bm{\alpha}}=\sup\left\{\nu>0\ |\ \limsup\limits_{k\to\infty}q_{k+1}^{\nu}\langle q_k\bm{\alpha}\rangle <\infty\right\}.$$
In particular, for any $\omega<\hat{\omega}_{\bm{\alpha}}$ the sequence $\left\{q_{k+1}^{\omega}\langle q_k\bm{\alpha}\rangle\right\}_{q\in\N}$ is bounded.
\end{Lem}
\begin{proof}
Let us denote $\hat{\omega}'_{\bm{\alpha}}=\sup\left\{\nu>0\ |\ \limsup\limits_{k\to\infty}q_{k+1}^{\nu}\langle q_k\bm{\alpha}\rangle <\infty\right\}$. We would like to show $\hat{\omega}_{\bm{\alpha}}=\hat{\omega}'_{\bm{\alpha}}$ and we note that the in particular part follows immediately from this equality. Suppose $\nu<\hat{\omega}_{\bm{\alpha}}$ then by the definition of $\hat{\omega}_{\bm{\alpha}}$, the system of inequalities \eqref{equ:dioexp} has nontrivial integer solutions in $q$ for all $t>0$ sufficiently large. In particular the system of inequalities $\langle q\bm{\alpha}\rangle <q^{-\nu}_{k+1}$, $|q|<q_{k+1}$ has nontrivial integer solutions in $q$ for all sufficiently large $k$. Since $|q|<q_{k+1}$ we have $\langle q_k\bm{\alpha}\rangle\leq \langle q\bm{\alpha}\rangle < q_{k+1}^{-\nu}$ for all sufficiently large $k$, implying that $\limsup\limits_{k\to\infty}q_{k+1}^{\nu}\langle q_k\bm{\alpha}\rangle<\infty$. Hence $\nu\leq \hat{\omega}'_{\bm{\alpha}}$ for all $\nu<\hat{\omega}_{\bm{\alpha}}$ and we have $\hat{\omega}_{\bm{\alpha}}\leq \hat{\omega}'_{\bm{\alpha}}$. For the other inequality, suppose $0<\nu<\hat{\omega}'_{\bm{\alpha}}$, then we have $\limsup\limits_{k\to\infty}q_{k+1}^{\nu}\langle q_k\bm{\alpha}\rangle=0$ implying that there exists some integer $k_0$ such that for all $k\geq k_0$ we have $\langle q_k\bm{\alpha}\rangle< q^{-\nu}_{k+1}$. Note that for any $t>q_{k_0}$ there exists some $k\geq k_0$ such that $q_k<t\leq q_{k+1}$. For such $t$, $q_k$ is a nontrivial integer solution to \eqref{equ:dioexp} since $\langle q_k\bm{\alpha}\rangle <q^{-\nu}_{k+1}\leq t^{-\nu}$ and $0<q_k<t$. This implies that $\hat{\omega}'_{\bm{\alpha}}\leq \hat{\omega}_{\bm{\alpha}}$ finishing the proof.
\end{proof}
\subsection{Effective density for irrational shifts}
In this subsection we state and prove our main result of this section. Let $\bm{\alpha}\in\R^n$ be an irrational vector, and let $\{q_k\}_{k\in\N}$ be the sequence of best simultaneous Diophantine approximation denominators of $\bm{\alpha}$. To state our result we first define a new Diophantine exponent attached to $\bm{\alpha}$ using the sequence $\{q_k\}_{k\in\N}$. Namely, let
$$\nu_{\bm{\alpha}}:=\inf\left\{\nu> 0\ |\ \limsup\limits_{k\to\infty}\frac{q_{k+1}}{q_k^{n\nu}}<\infty\right\}.$$
We now give the set of irrational vectors that we can handle. Let $n\geq 5$ and we define the set
$$\mathcal{DI}:=\left\{\bm{\alpha}\in \R^n\setminus \Q^n\ |\ \nu_{\bm{\alpha}}<\infty\ \textrm{and}\  \hat{\omega}_{\bm{\alpha}}>\frac{2}{n-2}\right\}.$$
We can now state our main result regarding the effective density for values of generic inhomogeneous quadratic forms with a fixed irrational shift.
\begin{Thm}\label{thm:moregeneral}
Let $n\geq 5$ and let $\bm{\alpha}\in \mathcal {DI}$ with $\nu_{\bm{\alpha}}$ and $\hat{\omega}_{\bm{\alpha}}$ as above. Then for any $\kappa\in \left(0,\frac{(n-2)\hat{\omega}_{\bm{\alpha}}-2}{n\left(1+\nu_{\bm{\alpha}}+\nu_{\bm{\alpha}}\hat{\omega}_{\bm{\alpha}}\right)}\right)$, for any $\xi\in \R$, and for almost every non-degenerate indefinite quadratic form $Q$ in $n$ variables, the system of inequalities \eqref{equ:maininequ} has integer solutions for all sufficiently large $t$.
\end{Thm}
\begin{rmk}\label{rmk:setdi}
We note that \lemref{lem:badlyapproximable} implies that $\nu_{\bm{\alpha}}\leq \omega_{\bm{\alpha}}$. Thus $\nu_{\bm{\alpha}}<\infty$ is a slightly weaker condition than $\bm{\alpha}$ is Diophantine and hence it is a full measure condition and the set $\mathcal{DI}$ contains all the vectors considered in \thmref{thm:irrational}. In contrast, since $\frac{2}{n-2}>\frac{1}{n}$ the condition $\hat{\omega}_{\bm{\alpha}}>\frac{2}{n-2}$ is very restrictive and is only satisfied by a null set. In fact, it was shown by German \cite[Theorem 3]{German2012} that $\frac{1}{n}\leq \hat{\omega}_{\bm{\alpha}}\leq 1$ for any irrational $\bm{\alpha}\in\R^n$ (see also \cite[Theorem 1.4]{Marnat2018}).  Hence when $n=3$ or $4$ the condition $\hat{\omega}_{\bm{\alpha}}>\frac{2}{n-2}$ is void and we note that this is exactly why our result can only handle the case when $n\geq 5$. On the other hand, when $n\geq 5$ let $l$ be a positive integer such that $l<\frac{n-2}{2}$. Let $\bm{\lambda}\in \R^l$ be any Diophantine vector and let $\textrm{Span}_{\Q}(\bm{\lambda})$ be the rational span of the coordinates of $\bm{\lambda}$ together with $1$. Then for any irrational $\bm{\alpha}\in \left(\textrm{Span}_{\Q}(\bm{\lambda})\right)^n$ such that its coordinates together with $1$ also span $\textrm{Span}_{\Q}(\bm{\lambda})$ we have $\nu_{\bm{\alpha}}\leq \omega_{\bm{\alpha}}=\omega_{\bm{\lambda}}<\infty$ and $\hat{\omega}_{\bm{\alpha}}= \hat{\omega}_{\bm{\lambda}}\geq \frac{1}{l}>\frac{2}{n-2}$, implying that $\bm{\alpha}\in \mathcal{DI}$. In addition to these examples which essentially come from lower dimensional Euclidean spaces, there are many other vectors lying in $\mathcal{DI}$. For instance Marnat and Moshchevitin  \cite[Theorem 1]{MarnatMoshchevitin2018} showed that for any $\frac{1}{n}\leq \omega<1$ there exist continuum many totally irrational Diophantine vectors $\bm{\alpha}\in \R^n$ such that $\hat{\omega}_{\bm{\alpha}}=\omega$. In particular this implies that there are continuum many totally irrational vectors lying in $\mathcal{DI}$.
\end{rmk}
\begin{proof}[Proof of \thmref{thm:moregeneral}]
We first prove this theorem for generic unit determinant forms. Fix an arbitrary compact set $\cK\subset G_n$ and fix integers $p_1,p_2$ with $p_1, p_2\geq 1$ and $p_1+p_2=n$. Let $Q_0$ be the fixed unit determinant form of signature $(p_1,p_2)$ as in \secref{sec:rational}. It suffices to prove the conclusion of \thmref{thm:moregeneral} for $g\cdot Q_0$ for $\mu_n$-almost every $g\in \cK$. For any $\xi\in \R$ and any $\kappa\in \left(0,\frac{(n-2)\hat{\omega}_{\bm{\alpha}}-2}{n\left(1+\nu_{\bm{\alpha}}+\nu_{\bm{\alpha}}\hat{\omega}_{\bm{\alpha}}\right)}\right)$, denote by $I_t=(\xi-t^{-\kappa},\xi+t^{-\kappa})$ and for any $g\in \cK$ let 
$$A_{g,I_t,t}=Q_0^{-1}(I_t)\medcap (B_t+\bm{\alpha})g$$ 
be as defined in \secref{sec:rational}. Let $\{q_k\}_{k\in\N}$ be the sequence of best simultaneous Diophantine approximation denominators of $\bm{\alpha}$ and let $\{\bm{r}_k\}_{k\in \N}$ be the sequence of partial convergents of $\bm{\alpha}$. Since $\kappa\in \left(0,\frac{(n-2)\hat{\omega}_{\bm{\alpha}}-2}{n\left(1+\nu_{\bm{\alpha}}+\nu_{\bm{\alpha}}\hat{\omega}_{\bm{\alpha}}\right)}\right)$, we can take $\nu>\nu_{\bm{\alpha}}$ sufficiently small and $\omega\in \left(\frac{2}{n-2},\hat{\omega}_{\bm{\alpha}}\right)$ sufficiently large such that $\kappa\in \left(0,\frac{(n-2)\omega-2}{n\left(1+\nu+\nu\omega\right)}\right)$. We fix these $\nu$ and $\omega$ for the rest of the proof and note that by the definition of $\nu_{\bm{\alpha}}$ and \lemref{lem:bestapp}, we have the sequences $\left\{\frac{q_{k+1}}{q_k^{n\nu}}\right\}_{k\in \N}$ and $\left\{q_{k+1}^{\omega}\langle q_k\bm{\alpha}\rangle\right\}_{k\in\N}$ are bounded (with the bounding constants depending on $\nu$, $\omega$ and $\bm{\alpha}$). In particular, we have
\begin{equation}\label{equ:di}
\|\bm{\alpha}-\bm{r}_k\|_{\infty}\ll_{\bm{\alpha},\omega} q_k^{-1}q_{k+1}^{-\omega}\ \textrm{for all $k\in \N$}.
\end{equation}For any $t>0$ let 
$$\cB_t:=\left\{g\in \cK\ |\ (\Z^n+\bm{\alpha})g\medcap A_{g,I_t,t}=\emptyset \right\}.$$
Then it suffices to show that $\limsup_{t\to \infty}\cB_t$ is of zero measure. 
Take $\{t_k=q_k^{\beta}\}_{k\in\N}$ for some $\beta\in \left(0,\frac{1+\omega}{1+\kappa}\right)$ to be determined. Since $\{t_k\}_{k\in\N}$ is unbounded we have
$$\limsup_{t\to\infty}\cB_t=\bigcap_{T>0}\bigcup_{t\geq T}\cB_t=\bigcap_{m\in \N}\bigcup_{k\geq m}\bigcup_{t_k\leq t< t_{k+1}}\cB_t.$$
Thus it suffices show the series $\sum_k\mu_n\left(\bigcup_{t_k\leq t<t_{k+1}}\cB_t\right)$ is summable for some choice of $\beta$. Suppose $g\in \bigcup_{t_k\leq t<t_{k+1}}\cB_t$, then there exists some $t_k\leq t<t_{k+1}$ such that $(\Z^n+\bm{\alpha})g\medcap A_{g,I_t,t}=\emptyset$. Since $t_k\leq t<t_{k+1}$ we have $A_{g,I_{t_{k+1}},t_k}\subset A_{g,I_t,t}$ implying that $(\Z^n+\bm{\alpha})g\medcap A_{g,I_{t_{k+1}},t_k}=\emptyset$. Shifting the set by $(\bm{r}_k-\bm{\alpha})g$ we get 
\begin{equation}\label{equ:relationinc}
(\Z^n+\bm{r}_k)g\bigcap \left(\left(Q_0^{-1}(I_{t_{k+1}})+(\bm{r}_k-\bm{\alpha})g\right)\medcap \left(B_{t_k}+\bm{r}_k\right)g\right)=\emptyset.
\end{equation}
Next, let $I_t'=(\xi-\frac{t^{-\kappa}}{2},\xi+\frac{t^{-\kappa}}{2})$, and we would like to show that there exists some $k_0\geq 1$ sufficiently large such that for any $k\geq k_0$ and for any $g\in \cK$ we have 
\begin{equation}\label{equ:relationir}
Q_0^{-1}(I_{t_{k+1}}')\medcap \left(B_{t_k}+\bm{r}_k\right)g\subset  \left(Q_0^{-1}(I_{t_{k+1}})+(\bm{r}_k-\bm{\alpha})g\right)\medcap \left(B_{t_k}+\bm{r}_k\right)g.
\end{equation}
We note that for any $\bm{v},\bm{w}\in \R^n$, $|Q_0(\bm{v})-Q_0(\bm{w})|\ll_n \|\bm{v}-\bm{w}\|_{\infty}\|\bm{v}+\bm{w}\|_{\infty}$. Hence for any $\bm{v}\in Q_0^{-1}(I_{t_{k+1}}')\medcap \left(B_{t_k}+\bm{r}_k\right)g$ with $\bm{w}:=\bm{v}-(\bm{r}_k-\bm{\alpha})g$ we have
$$\bigg||Q_0(\bm{v})-\xi|-|Q_0(\bm{w})-\xi|\bigg|\leq |Q_0(\bm{v})-Q_0(\bm{w})|\ll_n \|(\bm{r}_k-\bm{\alpha})g\|_{\infty}\|\bm{v}+\bm{w}\|_{\infty}\ll_{\cK}\|\bm{r}_k-\bm{\alpha}\|_{\infty}t_k,$$
where for the last inequality we used the assumption that $\bm{v}\in (B_{t_k}+\bm{r}_k)g$. Thus by \eqref{equ:di} we have for any $k\geq 1$ and for any $g\in \cK$
\begin{equation}\label{equ:intermed}
\bigg||Q_0(\bm{v})-\xi|-|Q_0(\bm{w})-\xi|\bigg|\ll_{n,\cK,\bm{\alpha},\omega}t_kq^{-1}_kq_{k+1}^{-\omega}= q_k^{\beta-1}q_{k+1}^{-\omega},
\end{equation}
where for the equality we used the assumption that $t_k=q_k^{\beta}$. 
Note that 
$$\frac{q_k^{\beta-1}q_{k+1}^{-\omega}}{t_{k+1}^{-\kappa}}=q_k^{\beta-1}q_{k+1}^{\beta\kappa-\omega}<q_k^{\beta-1-\omega+\beta\kappa}= o(1),$$
where for the inequality we used that $q_k<q_{k+1}$ and $\beta\kappa<\omega$ (since $\beta\kappa<\frac{1+\omega}{1+\kappa}\kappa<\omega$ with the second estimate following from $\kappa<\frac{(n-2)\omega-2}{n\left(1+\nu+\nu\omega\right)}<\omega$), and for the last estimate we used that $\beta-1-\omega+\beta\kappa< 0$ which follows from $\beta<\frac{1+\omega}{1+\kappa}$. This estimate, together with \eqref{equ:intermed} implies that there exists some $k_0\geq 1$ sufficiently large such that for any $k\geq k_0$ and for any $g\in \cK$ with $\bm{v}$ and $\bm{w}$ as above we have $\bigg||Q_0(\bm{v})-\xi|-|Q_0(\bm{w})-\xi|\bigg|<\frac{t_{k+1}^{-\kappa}}{2}$, implying \eqref{equ:relationir}. 

Next, let $\cO\subset G_n$ be an open neighbourhood of the identity element in $G_n$ satisfying $B_{\frac12t}\subset B_t h$ for any $h\in \cO$ and for any $t>0$ (for the existence of such an open neighborhood, cf. \cite[Section 2.2]{KelmerYu2018b}). Since $\cK$ is compact, there is a finite set $I\subset\cK$ such that $\cK\subset \bigcup_{x\in I}\cO x$. Then for any $g\in \cK$ there exists some $x\in I$ and $h\in \cO$ such that $g=hx$, implying that for $k\geq k_0$
\begin{equation}\label{equ:relationir2}
B_{\frac12(t_k-2\|\bm{\alpha}\|)}x\subset B_{t_k-2\|\bm{\alpha}\|}hx \subset (B_{t_k}+\bm{r}_k)g,
\end{equation}
where for the second inclusion we used that $B_{t_k-2\|\bm{\alpha}\|}\subset B_{t_k}+\bm{r}_k$ which follows from the triangle inequality and the estimate that $\|\bm{r}_k\|\leq 2\|\bm{\alpha}\|$. Now for any $x\in I$ and for any $k\geq k_0$ we denote by 
$$\underline{A}_{k,x}:=Q_0^{-1}(I_{t_{k+1}}')\medcap B_{\frac12(t_k-2\|\bm{\alpha}\|)}x,$$
and \eqref{equ:relationir} and \eqref{equ:relationir2} imply that for any $g\in \cK$, there exists some $x\in I$ such that $\underline{A}_{k,x}\subset\left(Q_0^{-1}(I_{t_{k+1}})+(\bm{r}_k-\bm{\alpha})g\right)\medcap \left(B_{t_k}+\bm{r}_k\right)g$. Now for any $k\geq k_0$ suppose $g\in \bigcup_{t_k\leq t<t_{k+1}}\cB_t$, together with \eqref{equ:relationinc} we get that there exists some $x\in I$ such that $(\Z^n+\bm{r}_k)g\medcap \underline{A}_{k,x}=\emptyset$. This implies that for $k\geq k_0$
$$\bigcup_{t_k\leq t<t_{k+1}}\cB_t\subset \bigcup_{x\in I}\left\{g\in \cK\ \left|\ (\Z^n+\bm{r}_k)g\medcap \underline{A}_{k,x}=\emptyset\right.\right\}\subset\bigcup_{x\in I}\cM^{(\cK,\bm{r}_k)}_{\underline{A}_{k,x},T_{k,x}}$$
where $T_{k,x}=\vol\left(\underline{A}_{k,x}\right)$. Thus applying \lemref{lemma:MAT} to $\cM^{(\cK,\bm{r}_k)}_{\underline{A}_{k,x},T_{k,x}}$ and applying \thmref{thm:volumeestimate} to $T_{k,x}=\vol\left(\underline{A}_{k,x}\right)$, we can bound for $k\geq k_0$
$$\mu_n(\bigcup_{t_k\leq t<t_{k+1}}\cB_t)\ll_{\cK,\bm{\alpha}}\sum_{x\in I}\frac{q_k^n}{T_{k,x}}\ll_{\cK}\frac{q_k^n}{t_{k+1}^{-\kappa}t_k^{n-2}}\ll_{\bm{\alpha},\nu} \frac{1}{q_k^{\beta(n-2)-n-n\nu\kappa \beta}},$$
where for the last estimate we used that $\{\frac{q_{k+1}}{q_k^{n\nu}}\}_{k\in \N}$ is bounded and that $t_k=q_k^{\beta}$. Since $q_k\geq c^k$ for some $c>1$, it suffices to show that there exists $\beta\in \left(0,\frac{1+\omega}{1+\kappa}\right)$ such that $\beta(n-2)-n-n\nu\kappa \beta>0$, or equivalently, $\beta>\frac{n}{n-2-n\nu\kappa}$. We note that since $\kappa\in \left(0,\frac{(n-2)\omega-2}{n\left(1+\nu+\nu\omega\right)}\right)$ we have $\frac{1+\omega}{1+\kappa}>\frac{n}{n-2-n\nu\kappa}$. Thus we can always take $\beta\in \left(0,\frac{1+\omega}{1+\kappa}\right)$ sufficiently large such that $\beta>\frac{n}{n-2-n\nu\kappa}$ finishing the proof for unit determinant forms. 

Finally for general determinants, for any $\lambda>0$ we want to show that for almost every $Q\in \cQ_{p_1,p_2}$ the system of inequalities 
\begin{equation}\label{equ:scaleineq}
|\lambda Q_{\bm{\alpha}}(\bm{v})-\xi|<t^{-\kappa}\quad  \textrm{and}\quad \|\bm{v}\|<t
\end{equation}
has integer solutions for all $t$ sufficiently large. Take some $\kappa'\in \left(\kappa,\frac{(n-2)\hat{\omega}_{\bm{\alpha}}-2}{n\left(1+\nu_{\bm{\alpha}}+\nu_{\bm{\alpha}}\hat{\omega}_{\bm{\alpha}}\right)}\right)$ and apply the above results for unit determinant forms for $\lambda^{-1}\xi$ and $\kappa'$ to get that for almost every $Q\in \cQ_{p_1,p_2}$ the system of inequalities 
$$|\lambda Q_{\bm{\alpha}}(\bm{v})-\xi|<\lambda t^{-\kappa'}\quad \textrm{and}\quad \|\bm{v}\|<t$$
has integer solutions for all $t$ sufficiently large. Since $\kappa'>\kappa$ we have $\lambda t^{-\kappa'}< t^{-\kappa}$ for all $t$ sufficiently large, implying that for almost every $Q\in \cQ_{p_1,p_2}$ the system of inequalities \eqref{equ:scaleineq} has integer solutions for all $t$ sufficiently large. This finishes the proof.
\end{proof} 

\begin{proof}[Proof of \thmref{thm:irrational}]
Let $\bm{\alpha}\in\R^n$ be as in this theorem, that is $\omega_{\bm{\alpha}}<\infty$ and $\hat{\omega}_{\bm{\alpha}}>\frac{2}{n-2}$. As mentioned in \rmkref{rmk:setdi} since $\nu_{\bm{\alpha}}\leq \omega_{\bm{\alpha}}$ we have $\bm{\alpha}\in\mathcal {DI}$. Then the theorem follows immediately from \thmref{thm:moregeneral} and noting that $\frac{(n-2)\hat{\omega}_{\bm{\alpha}}-2}{n\left(1+\omega_{\bm{\alpha}}+\omega_{\bm{\alpha}}\hat{\omega}_{\bm{\alpha}}\right)}\leq \frac{(n-2)\hat{\omega}_{\bm{\alpha}}-2}{n\left(1+\nu_{\bm{\alpha}}+\nu_{\bm{\alpha}}\hat{\omega}_{\bm{\alpha}}\right)}$.
\end{proof}


\end{document}